\documentstyle[12pt,amssymb,amscd]{amsart}  

\newcommand{\bbC}{{\Bbb C}}

\newcommand{\bbR}{{\Bbb R}}

\newcommand{\bbZ}{{\Bbb Z}}
\newcommand{\cF}{{\cal F}}

\newcommand{\cD}{{\cal D}}

\newcommand{\Sc}{{\cal S}}

\newcommand{\NC}{\operatorname{NC}}
\newcommand{\Mat}{\operatorname{Mat}}

\newcommand{\Ker}{\operatorname{Ker}}

\newcommand{\Tr}{\operatorname{Tr}}

\newcommand{\Ito}{\operatorname{Ito}}
\newcommand{\Str}{\operatorname{Str}}

\newcommand{\ad}{\operatorname{ad}}

\newcommand{\Lie}{\operatorname{Lie}}

\newcommand{\Herm}{\operatorname{Herm}}

\newcommand{\Gg}{{\frak g}}
\newcommand{\GF}{{\frak F}}

\newcommand{\Hol}{\operatorname{Hol}}

\title{Noncommutative geometry and path integrals}

\author{M. Kapranov}
\email{mikhail.kapranov@@yale.edu}
\address{Department of Mathematics, Yale University, 10 Hillhouse Avenue,
New Haven, CT 06520}
\begin{document}
 \maketitle
 
  \begin{center}
\vspace*{.7cm}
{\it To Yuri Ivanovich Manin on his 70th birthday}
\end{center}

\vskip 1cm

\centerline {\bf Introduction}

\vskip 1cm

\noindent {\bf (0.1)}   A monomial in noncommutative variables $X$ and $Y$, say, 
$X^iY^jX^kY^l...$, 
 can be visualized as a lattice path in the plane, starting from $0$, 
 going $i$ steps
in the horizontal direction, $j$ steps in the vertical one , then again $k$ steps in
the horizontal one, and so on. Usual commutative monomials are often visualized as lattice
points, for example $x^a y^b$ corresponds to the point $(a,b)$. To lift such a monomial
to the noncommutative domain, is therefore the same as to choose a ``history" for $(a,b)$,
i.e., a lattice path originating at 0 and ending at $(a,b)$.

\vskip .2cm

This correspondence between paths and noncommutative monomials can be extended to more
general piecewise smooth paths, if we deal with exponential functions instead.
Let us represent our commutative variables as $x=e^z, y=e^w$, then a monomial
will be replaced by the exponential $e^{az+bw}$ and we are free to take  $a$ and $b$ to be any real 
numbers. To lift this exponential to the noncommutative domain, i.e. to a series
in $Z,W$ where $X=e^Z, Y=e^W$, one needs to choose a path $\gamma$ in $\bbR^2$ joining 0 with $(a,b)$.
One can easily see this by approximating $\gamma$ by lattice paths with step
 $1/M$, $M\to\infty$,
and working with monomials in $X^{1/M}= e^{Z/M}$ and $Y^{1/M}=e^{W/M}$. Denote this exponential series
$E_\gamma(Z,W)$. 

\vskip .2cm

This suggests the possibility of a ``noncommutative Fourier transform" (NCFT) identifying appopriate spaces of
functions of noncommuting variables (say, of matrices of intederminate size) with spaces of ordinary
functions or measures on the space of paths. For example, to a measure $\mu$ on the space $\Pi$ of paths
(or some completion of it) we want to associate the function of  $Z,W$ 
given by
$$\cF(\mu)(Z,W) = \int_{\gamma\in\Pi} E_\gamma(Z,W) d\mu(\gamma),
\leqno (0.1.1)$$
The basic phenomenon here seems to be that the
two types of functional spaces (noncommutative functions of $n$ variables vs. ordinary commutative functions
 but on the space paths in $\bbR^n$), have, on some fundamental level,
{\em the same size}.

The goal of this paper and the ones to follow [K1-2] is to investigate this idea from several points of view.

\vskip .3cm

\noindent {\bf (0.2)} The concept  of NCFT seems to 
implicitly underlie
the very foundations of quantum mechanics such as  the equivalence of the Lagrangian and
Hamiltonian approaches to the theory. Indeed, the Lagrangian point of view deals with path integrals
while the Hamiltonian one works with noncommuting operators.  
Further, it is very close to the concept of the ``Wilson loop" functional
(trace of the holonomy)  in Yang-Mills theory [Po].
Note that the exponential $E_\gamma$, being itself the holonomy of a certain
formal connection, is invariant under reparametrization of the path.
Quantities invariant under reparametrization are particularly important in string theory,
and the reparametrization invariance of the Wilson loop led to conjectural relations
between strings and $N\to\infty$ limit of Yang-Mills  theory [Po].  

As the integral transform $\cF$ should, intuitively, act between spaces of the same size,
it does not lead to any loss of information and can therefore be viewed as ``path integration 
without integration". The actual integration occurs when we restrict the function $\cF(\mu)$
to the commutative locus, i.e., make $Z$ and $W$ commute. Alternatively, 
instead of allowing $Z,W$ to be arbitrary matrices, we take them to be scalars. 
Then all paths having the same endpoint
will contribute to make up a single  Fourier mode of the commutativized function.
We arrive at the following conclusion: {\em the natural homomorphism
$R\to R_{ab}$ of a noncommutative ring to its maximal commutative quotient
is the algebraic analog of path integration}. 

\vskip .3cm

\noindent {\bf (0.3)} The idea that the space of paths is related to the free group and to
its various versions has been clearly enunciated by K.-T. Chen [C1] in the 1950's and can be traced
throughout almost all of his work [C0]. Apparently, 
 a lot more can be said about this classical subject.  Thus, the universal connection with
 values in the free Lie algebra (known to Chen and appearing in (2.1) below) leads to  beautiful non-holonomic
 geometry on the free nilpotent Lie groups $G_{n,d}$, which is still far from being fully understood,
 see [G]. 
 
 Well known examples of measures on path spaces are provided by probability theory and we spend some time
 in \S 4 below to formulate various results from probabilistic literature in terms of NCFT.
 Most importantly, the Fourier transform of the Wiener measure on paths  in
 $\bbR^n$ is the noncommutative Gaussian series
 $\exp\left( -\sum_{i=1}^n Z_j^2\right)$ where $Z_i$ are considered as noncommuting variables. 
 We should mention here the recent book by Baudoin [Ba] who considered the idea of associating
 a noncommutative series to a stochastic process. 
 It is clearly the same type of construction as our NCFT except in the
  framework of probability theory:
 parametrized paths, positive measures etc.

\vskip .3cm

\noindent {\bf (0.4)} I would like to thank R. Beals, E. Getzler, H. Koch, Y.I. Manin and M. A. Olshanetsky
for useful discussions. This paper was written during my stay at the Max-Planck-Institute f\"ur Mathematik 
in Bonn and I am grateful to the Institute for support and excellent working conditions.
This work was also partially supported by an NSF grant. 

\vfill\eject

\centerline{\bf 1. Noncommutative monomials and lattice paths.}

\vskip 1cm

\noindent {\bf (1.1) Noncommutative polynomials and the free semigroup.}
Consider 
$n$ noncommuting (free) variables $X_1, ..., X_n$ and form the algebra 
of noncommutative polynomials in these variables. 
This algebra will be denoted by 
${\bbC}\langle X_1, ..., X_n\rangle$. It is the same as the tensor algebra
$$T(V) = \bigoplus_{d=1}^\infty V^{\otimes d}, \quad V = \bbC^n = \bigoplus \bbC\cdot X_i.$$
A noncommutative monomial in $X= (X_1, ..., X_n)$ is, as described in the
Introduction, the same as a monotone lattice path in ${\bbR}^n$ starting 
at 0. We denote by $F^+_n$ the set of all such paths 
and write $X^\gamma$ for the monomial corresponding to a path $\gamma$. 
The set $F_n^+$ is a semigroup  with the following operation. If $\gamma,
\gamma'$ are two monotone paths as above starting at 0, then
$\gamma\circ\gamma'$ is obtained by translating $\gamma$ so that its
beginning meets the end of $\gamma'$ and then forming the
composite path. It is clear that $F_n^+$ is 
the free semigroup on $n$ generators. 
Thus a typical noncommutative polynomial is written as 
$$f(X_1, ... X_n) = f(X) = \sum_{\gamma\in F^+_n} a_\gamma X^\gamma.
\leqno (1.1.1) $$ 
Along with ${\bbC}\langle X_1, ..., X_n\rangle$ we will consider the algebra 
${\bbC}[x_1, ..., x_n]$ of usual (commutative) polynomials in the variables 
$x_1, ..., x_n$. A typical such polynomial will be written as 
$$g(x_1, ..., x_n) = g(x) = \sum_{\alpha\in {\bf Z}_+^n} b_\alpha x^\alpha, 
\quad  x^\alpha = x_1^{\alpha_1}... x_n^{\alpha_n}.\leqno (1.1.2)$$ 
The two algebras are related by the 
 {\it commutativization 
homomorphism} 
$$c: {\bbC}\langle X_1, ..., X_n\rangle \to {\bbC}[x_1, ..., x_n],
\leqno (1.1.3)$$ 
which takes $X_i$ to $x_i$. 
For a path $\gamma\in\Gamma_n$ let $e(\gamma)\in {\bbZ}_+^n$ denote the 
end point of $\gamma$. Then we have 
$$c(X^\gamma) = x^{e(\gamma)}.\leqno (1.1.4)$$ 
This means that at the level of coefficients the 
commutativization homomorphism is given by the summation over paths 
with given endpoints: if $g(x) = c(f(X))$, then 
$$b_\alpha = \sum_{e(\gamma)=\alpha} a_\gamma.\leqno (1.1.5)$$ 

\vskip .3cm

\noindent {\bf (1.2) Noncommutative power series.} 
  Let 
 $I\subset {\bbC}\langle X_1, ..., X_n\rangle$ 
be the span of monomials of degree $\geq 1$. Then clearly $I$ is a 2-sided 
ideal and $I^d$ is the span of monomials of degree $\geq d$. We define the 
algebra ${\bbC}\langle\langle X_1, ..., X_n\rangle\rangle$ as 
 the completion of 
${\bbC}\langle X_1, ..., X_n\rangle$ in the $I$-adic topology. Explicitly, 
elements of ${\bbC}\langle\langle X_1, ..., X_n\rangle\rangle$ can be seen 
as infinite formal linear combinations of noncommutative monomials, i.e.,
expressions of the form $\sum_{\gamma\in F_N^+} a_\gamma X^\gamma$. For example,
$$e^{X_1}\cdot e^{X_2} = \sum_{i,j=0}^\infty {X_1^i X_2^j\over i! j!}, 
 \quad  {1\over 1-(X_1+X_2)} = \sum_{\gamma\in F^+_2} X^\gamma
 \leqno (1.2.1)$$
 are noncommutative power series. We will be also interested in convergence of noncommutative
 series. A series
 $f(X) =\sum_{\gamma\in F_n^+} a_\gamma X^\gamma$ will be called {\em entire},
 if
 $$\lim_{\gamma\to\infty} R^{l(\gamma)} |a_\gamma| = 0, \quad \forall R>0. \leqno (1.2.2)$$ 
 Here $l(\gamma)$ is the length of the path $\gamma$ and the limit is taken over
 the countable set $F_n^+$ (so no ordering of this set is needed). If (1.2.2) is the case,
 then for any $N$ and for any square matrices $X_1^0, ..., X_n^0$ of size $N$ the
 series of matrices $\sum a_\gamma (X^0)^\gamma$ obtained by specializing $X_i\to X_i^0$,
 converges absolutely. We denote by 
 $\bbC\langle\langle X_1, ..., X_n\rangle\rangle^{ent}$ the set of entire series.
 It is clear that this set is a subring. 
 
 \vskip .3cm
 
 \noindent {\bf (1.3) Noncommutative Laurent polynomials.} By a noncommutative 
Laurent monomial in $X_1, ..., X_n$ we will mean a monomial in positive and 
negative powers of the $X_i$ such as, e.g., $X_1X_2 X_1^{-1} X_2^5$. 
In other words, this is an element of $F_n$, the free noncommutative 
group on the generators $X_i$. A noncommutative Laurent polynomial 
is then a finite formal linear combination of such monomials 
i.e., an element of the group algebra of $F_n$. We will denote this 
algebra by 
$${\bbC}\langle X_1^{\pm 1}, ..., X_n^{\pm 1}\rangle = {\bbC}[F_n].\leqno
(1.3.1)$$ 
As before, a noncommutative Laurent monomial corresponds to a lattice 
path in ${\bbR}^n$ beginning at 0 but not necessarily monotone. 
These paths are defined up to cancellation of pieces consisting 
of a sub-path and the same  sub-path run in the opposite direction 
immediately afterwards. 

We retain the notation $X^\gamma$ for the monomial corresponding to 
a path $\gamma$. We also write $(-\gamma)$ for the path inverse to
$\gamma$, so $X^{-\gamma} = (X^\gamma)^{-1}$. 

\vskip .3cm

\noindent {\bf (1.4) Noncommutative Fourier transform: discrete case.} 
The usual 
(commutative) Fourier transform relates the spaces of functions on a locally 
compact abelian group $G$ and its Pontryagin dual $\widehat{G}$. The  
``discrete" case
 $G= \bbZ^n$, $\widehat{G} = {(S^1)}^n$ corresponds to the theory of 
Fourier series.   

In the algebraic formulation, the discrete Fourier transform identifies the 
space of finitely supported functions 
$$b: {\bbZ}^n \to {\bbC}, \quad \alpha\mapsto b_\alpha, \quad |{\rm Supp}(b)| 
< \infty, \leqno (1.4.1)$$ 
with the space ${\bbC}[x_1^{\pm 1}, ..., x_n^{\pm 1}]$ of Laurent polynomials. 
It is given by the well known formulas 
$$(b_\alpha) \mapsto f, \quad f(x) = \sum_{\alpha\in {\bbZ}^n} b_\alpha 
x^\alpha,\leqno (1.4.2) $$ 
$$f\mapsto (b_\alpha), \quad b_\alpha = \int_{|x_1| =... = |x_n|=1} 
f(x) x^{-\alpha} d^*x_1 ... d^*x_n,\leqno (1.4.3) $$ 
where $d^*x$ is the Haar measure on $S^1$ with volume 1. 
Our goal in this section is to give a generalization of these formulas for 
noncommutative Laurent polynomials. 

Instead of (1.4.1) we consider the space of finitely supported functions
$$a: F_n\to \bbC, \quad \gamma\mapsto a_\gamma, 
\quad |{\rm Supp}(a)|<\infty.  \leqno (1.4.4)$$
The discrete noncommutative Fourier transform is the identification of this
space with $\bbC \langle X_1^{\pm 1}, ..., X_n^{\pm 1}\rangle$ via
 $$(a_\gamma)\mapsto f, \quad f(X) = \sum_{\gamma\in F_n} a_\gamma X^\gamma.
 \leqno (1.4.5) $$ 
This identification ceases to look like a tautology if we regard
a noncommutative Laurent polynomial as a function $f$  which to any 
$n$ invertible
elements $X_1^0, ..., X_n^0$ of any associative algebra $A$ 
associates an element $f(X_1^0, ..., X_n^0)\in A$. We want then to recover
the coefficients $a_\gamma$ in terms of the values of $f$ on various elements
of various $A$. Most importantly, we will consider $A= {\rm Mat}_N(\bbC)$,
the algebra of matrices of size $N$ and let $N$ be arbitrary. To get a generalization of
(1.4.3) we replace the unit circle $|x|=1$ by the group of unitary matrices
$U(N)\subset {\rm Mat}_N(\bbC)$. 
 Let $d^*X$ be the Haar measure on 
$U(N)$ of volume 1.

 The following result is a consequence
of the so-called  ``asymptotic freedom theorem for unitary matrices" 
due to Voiculescu [V], see also [HP] for a more elementary exposition. 

\proclaim (1.4.6) Theorem.
 If $f(X) = \sum_{\gamma\in F_n} a_\gamma X^\gamma$ is a 
noncommutative 
Laurent polynomial, then we have 
$$a_\gamma = \lim_{N\to\infty} {1\over N}  \,\, {\rm tr} \int_{X_1, ..., X_n\in U(N)} 
f(X_1, ..., X_n)\,  X^{-\gamma} \,\, d^*X_1 ... d^* X_n.$$

As for the commutative case, the theorem is equivalent to the following 
ortogonality 
relation. It is this relation that is usually called the 
"asymptotic freedom" in the literature. 

\proclaim (1.4.7) Reformulation. Let $\gamma\in F_n$ be a nontrivial 
  lattice path. Then 
$$\lim_{N\to\infty} {1\over N} \,\, {\rm tr} \int_{X_1, ..., X_n\in U(N)} 
X^\gamma \,\, d^*X_1 ... d^*X_n = 0.
 $$ 
Note that for $\gamma=0$ the integral is equal to 1 for any $N$. 

Passing to the $N\to\infty$ limit is unavoidable here since for any given $N$ there
exist noncommutative polynomials which vanish identically on ${\rm Mat}_N(\bbC)$.
An example is provided by the famous Amitsur-Levitsky polynomial
$$f(X_1, ..., X_{2N}) = \sum_{\sigma\in S_{2N}} {\rm sgn}(\sigma) X_{\sigma(1)} \cdot ... 
\cdot X_{\sigma(2N)}.$$


\vfill\eject

\centerline{\bf 2. Noncommutative exponential functions.}

\vskip 1cm

\noindent {\bf (2.1) The universal connection and noncommutative exponentials.}
Let us introduce the ``logarithmic" variables $Z_1, ..., Z_n$, so that we have the embedding
$$\bbC\langle X_1, ..., X_n\rangle \subset \bbC\langle\langle Z_1, ..., Z_n\rangle\rangle,
\quad X_i \mapsto e^{Z_i}. \leqno (2.1.1)$$
The algebra $\bbC\langle\langle Z_1, ..., Z_n\rangle\rangle$ is a projective limit of
finite-dimensional algebras, namely
$$\bbC\langle\langle Z_1, ..., Z_n\rangle\rangle\quad = \quad \lim_{\longleftarrow}{}_d\,\,\,
\bbC\langle Z_1, ..., Z_n\rangle/I^d, \leqno (2.1.2)$$
where the ideal $I$ is as in (1.2). 

Consider the space $\bbR^n$ with coordinates $y_1, ..., y_n$. On this space we have the
following 1-form with values in $\bbC\langle\langle Z_1, ..., Z_n\rangle\rangle$:
$$\Omega \quad = \quad \sum_i Z_i dy_i \quad \in\quad \Omega^1(\bbR^n) \otimes
\bbC\langle\langle Z_1, ..., Z_n\rangle\rangle. \leqno (2.1.3)$$
We consider the form as a connection on $\bbR^n$.  One can see it as the
universal translation invariant connection on $\bbR^n$, an algebraic version of
the connection of Kobayashi on the path space, see [Si], \S 3. 

Let $\gamma$ be any piecewise smooth path in $\bbR^n$. We define the 
noncommutative exponential function corresponding to $\gamma$ to be the
holonomy of the above connection along $\gamma$:
$$E_\gamma(Z) = E_\gamma(Z_1, ..., Z_n) = P \exp \int_\gamma \Omega \quad\in\quad \bbC\langle\langle Z_1, ...,
 Z_n\rangle\rangle.
\leqno (2.1.4)$$
The holonomy can be understood by passing to finite-dimensional quotients as in (2.1.2)
and solving an ordinary differential equation with values in each such quotient. 
It is clear that $E_\gamma (Z)$ becomes unchanged under parallel translations of $\gamma$,
since the form $\Omega$ is translation invariant. So in the following we will always assume
that $\gamma$ begins at $0$.

Further, $E_\gamma(Z)$ is invariant under reparametrizations of $\gamma$: this is a general
property of the holonomy of any connection. So let us give the following
definition. 

\proclaim (2.1.5) Definition. Let $M$ be a $C^\infty$-manifold.
An (oriented) unparametrized path in $M$ is an equivalence class of 
pairs $(I, \gamma: I\to M)$ where $I$ is a smooth manifold with boundary 
diffeomorhic to $[0,1]$ and $\gamma$ is a piecewise smooth map $I\to M$. 
Two such pairs $(I,\gamma)$ and $(I', \gamma')$ are equivalent if there 
is an orientation preserving piecewise smooth homeomorphism $\phi: I\to I'$ such that 
$\gamma=\gamma'\circ 
\phi$. 

We will denote an unparametrized path simply by $\gamma$. 

\vskip .2cm

\noindent {\bf (2.1.6) Example.} Let $\gamma$ be a straight segment in ${\bbR}^2$ 
joining (0,0) and (1,1). 
Let also $\delta$ be the path consisting of the horizontal segment $[(0,0), 
(0,1)]$ and the 
vertical segment $[(0,1), (1,1)]$. Let $\sigma$ be the path consisting of the vertical segment
$[(0,0), (1,0)]$ and the horizontal segment $[(1,0), (1.1)]$.  Then 
$$E_\gamma(Z_1, Z_2) = e^{Z_1+Z_2}, \quad E_\delta(Z_1, Z_2) = e^{Z_1} 
e^{Z_2}, \quad E_\sigma(Z) = e^{Z_2} e^{Z_1}.$$ 
More generally, if $\gamma$ is a lattice path corresponding to the
integer lattice $\bbZ^n$, then $E_\gamma(Z) = X^\gamma$ is the
noncommutative monomial in $X_i = e^{Z_i}$ associated to $\gamma$ as
in \S 1. 

\vskip .2cm

Let $\gamma, \gamma'$ be two unparametrized paths in $\bbR^n$ starting at 0.
Their product $\gamma\circ\gamma'$ is the path obtained by translating
$\gamma$ so that its beginning meets the end of $\gamma'$ and then
forming the composite path. The set of $\gamma$'s with this operation
forms a semigroup. For a path $\gamma$ we denote by $\gamma^{-1}$ the path
obtained by translating $\gamma$ so that its end meets $0$ and then
taking it with the opposite orientation. Finally, we denote by $\Pi_n$
the set of paths as above modulo cancellations, i.e., forgetting
sub-paths of a given path consisting of a segment and then immediately
of the same segment run in the opposite direction. Clearly the set $\Pi_n$
forms a group which we will call {\em the group of paths} in $\bbR^n$. 

The standard properties of the holonomy of connections imply the following:

\proclaim (2.1.7) Proposition. (a) We have 
$$E_{\gamma\circ\gamma'}(Z) = 
E_\gamma(Z)\cdot E_{\gamma'}(Z), \quad E_{\gamma^{-1}}(Z) = E_\gamma(Z)^{-1}$$
(equalities in  $\bbC\langle\langle Z_1, ..., Z_n\rangle\rangle$). 
\hfill\break
(b) The series $E_\gamma(Z)$ is entire, i.e., it converges for any given 
$N$ by $N$ matrices $Z_1^0, ..., Z_n^0$. \hfill\break
(c) If $Z_1^0, ..., Z_n^0$ are Hermitian, then $E_\gamma(iZ^0_1, ..., iZ^0_n)$
is unitary.

The property (a) implies that $E_\gamma(Z)$ depends only on the image of
$\gamma$ in the group $\Pi_n$.
Further, let us consider the commutativisation homomorphism
$$c: \bbC\langle\langle Z_1, ..., Z_n\rangle\rangle \to \bbC[[z_1, ..., z_n]].
\leqno (2.1.8)$$

The following is also obvious.

\proclaim (2.1.9) Proposition. If $a = (a_1, ..., a_n)$ is the endpoint of
$\gamma$, then
$$c(E_\gamma(Z)) = e^{(a,z)} $$
is the usual exponential function.

Thus there are as many ways to lift $e^{(a,z)}$ into the noncommutative
domain as there are paths in $\bbR^n$ joining $0$ and $a$.

\vskip .3cm

\noindent {\bf (2.2) Idea of a noncommutative Fourier transform.} 
The above observations suggest that there should be a version of Fourier transform
which would identify an appropriate space of measures on $\Pi_n$ with
an appropriate space of functions of $n$ noncommutative variables $Z_1, ..., Z_n$,
via the formula
$$\mu \mapsto f(Z_1, ..., Z_n) = \int_{\gamma\in \Pi_n} 
E_\gamma(iZ_1, ..., iZ_n) \cD\mu(\gamma).
\leqno (2.2.1)$$
The integral in (2.2.1) is thus a path integral. 
The concept of a ``function of noncommuting variables" is of course open to interpretation.
Several such interpretations are currently being considered in Noncommutative Geometry.

In the present paper we adopt a loose point of view that a function of $n$ noncommutative
variables is an element of an algebra $R$ equipped with a homomorphism $\bbC\langle Z_1, ..., 
Z_n\rangle\to R$. We will assume that this homomorphism realizes $R$ as some kind of
completion, or localization (or both) of $\bbC\langle Z_1, ..., Z_n\rangle$. In other
words, that $R$ does not have ``superfluous" elements, 
independent of the images of the $Z_i$. See [Ta] for an early attempt to define noncommutative
functions in the analytic context. 

\vskip .2cm

\noindent {\bf (2.2.2) Examples.} We can take $R=\bbC\langle \langle Z_1, ..., Z_n\rangle \rangle^{ent}$,
the algebra of entire power series. Alternatively we can take $R$ to be the skew field
of ``noncommutative rational functions" in $Z_1, ..., Z_n$ constructed by P. Cohn [Coh]. 
Thus expressions such as
$$\exp(Z_1^2 + Z_2^2), \quad \left(Z_1^2 + Z_2^2\right)^{-1}, \quad (Z_1Z_2-Z_2Z_1)^{-1}+
Z_3^{-2}Z_1$$
are considered as noncommutative functions. 

\vskip .2cm

It will be important for us to be able to view a ``function" $f(Z_1, ..., Z_n)$ as above
as an actual function defined on appropriate subsets of $n$-tuples of $N$ by $N$ matrices
for each $N$ and taking values in matrices of the same size. 

\vskip .2cm

Similarly, the group $\Pi_n$ can also possibly be replaced
 by various related objects (completions).  In  this paper we will consider
  several approaches such as  completion by  a pro-algebraic group or completion
  by continuous paths. 
  
  \vskip .2cm
 
 Alternatively, functions on $\Pi_n$ should correspond to ``noncommutative measures", 
 or distributions on the space of noncommutative functions. Examples of
 such ``measures" are being studied in the Free  Probability Theory 
 [HP] [NS] [VDN]. See \S 6 below. 
 
 Note that we have a surjective homomorphism of groups
 $$e: \Pi_n\to \bbR^n, \quad \gamma\mapsto e(\gamma). \leqno (2.2.3)$$
Here $e(\gamma)$ is the endpoint of $\gamma$. 
 One important property of the NCFT is the following principle which is just a consequence
 of Proposition 2.1.9: under the Fourier transform the integration over paths
 with given beginning and end, i.e., the pushdown of measures on $\Pi_n$ to
 measures on $\bbR^n$ corresponds to a simple algebraic operation: the commutativization
 homomorphism
 $$c: R \to R/([R,R]),\leqno (2.2.4)$$
 where $R$ is a noncommutative algebra and the RHS is the 
 maximal commutative quotient of $R$. 
 
 \vskip .3cm
 
 \noindent {\bf (2.3) Relation to Chen's iterated integrals.} 
 Let us recall the main points of Chen's theory.
 Let $M$ be a smooth 
manifold, $\gamma$ an unparametrized  path and $\omega$ a 1-form on $M$. 

Along with the "definite integral" $\int_\gamma\omega$, we can consider the 
"indefinite integral" which is a function ``on $\gamma$", or, more 
precisely, on the abstract interval $I$ such that $\gamma$ is a map 
$I\to M$.  For any $t\in I$ we have the sub-path $\gamma_{\leq t} $ 
going from the beginning of $I$ until $t$ and we have the function 
$$\int_{(\gamma)}\omega: \, I\to {\bbC}, \quad t\mapsto \int_{\gamma_{\leq 
t}}\omega.$$ 
If now $\omega_1$ and $\omega_2$ are two 1-forms on $M$, we can form a new 
1-form 
on $\gamma$ by multiplying (the restriction of) $\omega_2$ and the function 
$\int_{(\gamma)} \omega_1$. Then this form can be integrated along $\gamma$. 
The result is called the {\it iterated integral} 
$$\int_\gamma\biggl( \omega_2\cdot \int_{(\gamma)}\omega_1\biggl).$$ 
Note that if we think of $\gamma$ as a map $\gamma:I\to M$, then 
the iterated integral is equal to 
$$\int_{t_1\leq t_2\in I} \gamma^*(\omega_1)(t_1) \gamma^*(\omega_1)(t_2).$$ 
Note that integration over all $t_1, t_2\in I$ would give the product 
$\bigl(\int_\gamma \omega_1\bigr)\cdot\bigl(\int_\gamma\omega_2\bigr)$.

Similarly, one defines the $d$-fold iterated integral of $d$  1-forms $\omega_1, 
..., 
\omega_d$ on $M$ by induction: 
$$\int^\rightarrow_\gamma \omega_1\cdot ...\cdot \omega_d = \int_\gamma 
\biggl(\omega_d\cdot \int^\rightarrow_{(\gamma)} \omega_1\cdot ...\cdot 
\omega_{d-1} 
\biggr),$$ 
where the $(d-1)$-fold indefinite iterated integral is defined as the function 
on $I$ of the form 
$$t\to\int_{\gamma_{\leq t}}^\rightarrow \omega_1\cdot ...\cdot \omega_{d-1}.$$ 
As before the iterated integral is equal to the integral over the $d$-simplex: 
$$\int^\rightarrow_\gamma \omega_1\cdot ...\cdot \omega_d = 
\int_{t_1\leq ...\leq t_d\in I} \gamma^*\omega_1(t_1) 
...\gamma^*\omega_d(t_d).$$

\vskip .1cm 

The concept of iterated integrals extends in an obvious way to 1-forms with 
values in any associative (pro-)finite-dimensional $\bbC$-algebra $R$. The 
 well known Picard 
series for the holonomy of a connection consists  exactly  of such iterated 
integrals. 
We state this as follows. 

\proclaim (2.3.1) Proposition. Let $R$ be any (pro-)finite-dimensional 
associative $\bbC$-algebra, and 
 $\Omega$ be a 1-form on $M$ with values in $R$ 
considered as a connection form. Then the parallel transport along an 
unparametrized 
path $\gamma$ has the form 
$$P\exp\int_\gamma \Omega = \sum_{d=0}^\infty \int^\rightarrow_\gamma 
\Omega\cdot ... 
\cdot\Omega.$$ 

Let us specialize this to $M={\bbR}^n$, $R= {\bbC}\langle\langle Z_1, ..., Z_n 
\rangle\rangle$ and $\Omega = \sum Z_i dy_i$. We obtain: 

\proclaim (2.3.2) Corollary. The coefficient of the series $E_\gamma(Z_1, ..., Z_n)$ at 
any noncommutative monomial $Z_{i_1} ... Z_{i_d}$ is equal to the iterated 
integral 
$$\int^\rightarrow_\gamma dy_{i_1}\cdot ...\cdot dy_{i_d}.$$ 

Thus $E_\gamma$ is the generating function for all the  iterated integrals 
involving 
constant 1-forms on ${\bbR}^n$.

\vskip .2cm

\noindent {\bf (2.3.3) Example.} By the above 
$$E_\gamma(Z) = 1 + \sum a_i Z_i + \sum b_{ij} Z_iZ_j + ...$$ 
where $a_i = \int_\gamma dy_i$ is the $i$th coordinate of the endpoint of 
$\gamma$ and 
$$b_{ij} = \int_\gamma\biggl( dy_i \cdot\int_{(\gamma)} dy_j\biggr) 
=\int_\gamma y_j dy_i.$$ 
Suppose that $\gamma$ is closed, so $a_i=0$. Then $b_{ii}=0$ and 
for $i\neq j$ we have that $b_{ij}$ is the oriented area encirlced by $\gamma$ 
after the projection to the $(i,j)$-plane. 

\vskip .2cm

The following was proved by Chen [C2].

\proclaim (2.3.4) Theorem. The homomorphism $\Pi_n\to
\bbC \langle\langle Z_1, ..., Z_n\rangle\rangle ^*$ sending
$\gamma$ to $E_\gamma$ is injective. In other words, if a path $\gamma$
has all iterated integrals as above equal to 0, then $\gamma$ is
(equivalent modulo cancellations to) 
a constant path (situated at 0). 

\vskip .3cm

\noindent {\bf (2.4) Group-like and primitive elements.} Let
$FL(Z_1, ..., Z_n)$ be the free Lie algebra generated by $Z_1, ..., Z_n$.
It is characterized by the obvious universal property, see [R]
for background. This property implies that we have a Lie algebra
homomorphism
$$h: FL(Z_1, ..., Z_n) \to {\bbC}\langle Z_1, ..., Z_n\rangle,\leqno
(2.4.1)$$ 
and this homomorphism identifies ${\bbC}\langle Z_1, ..., Z_n\rangle$
with the universal enveloping algebra of $FL(Z_1, ..., Z_n)$. 
Further, let us consider the 
Hopf algebra structure on
${\bbC}\langle Z_1, ..., Z_n\rangle$
given on the generators by
$$\Delta(Z_i) = Z_i\otimes 1 + 1\otimes Z_i.\leqno (2.4.2)$$
The following result, originally due to K. Friedrichs, is a particular
case of a general property of enveloping algebras. 

\proclaim (2.4.3) Theorem.  
The image of $h$ consists precisely of all  primitive elements,
i.e., of elements $f$ such that 
$\Delta(f) = f\otimes 1 + 1\otimes f$. 

 We will also use the term {\em Lie elements} for primitive
 elements of $ \bbC\langle Z_1, ..., Z_n\rangle$.  

\vskip .2cm

Further, consider the noncommutative power series algebra
$\bbC\langle\langle Z_1, ..., Z_n\rangle\rangle$. It is
naturally a topological Hopf algebra with respect to the comultiplication
given  by (2.4.2) on generators and extended by additivity,
multiplicativity and continuity. 

The free Lie algebra is 
graded: 
$$FL(Z_1, ..., Z_n) = \bigoplus_{d\geq 1} FL(Z_1, ..., Z_n)_d,
\leqno (2.4.4)$$ 
where $FL(Z_1, ..., Z_n)_d$ is the span of Lie monomials containing exactly $d$ 
letters. 
We denote by 
$${\Gg}_n  = \prod_{d\geq 1} FL(Z_1, ..., Z_n)_d
\leqno (2.4.4)$$ its completion, i.e., 
the set of formal {\it Lie series}. This is a complete topological Lie algebra.
We clearly have an embedding of $\Gg_n$ into 
$\bbC\langle\langle Z_1, ...,
Z_n\rangle\rangle$ induced by the embedding of the graded components as above.
Further, degree-by-degree considerations and Theorem 2.4.3 imply the following:

\proclaim (2.4.5) Corollary. A noncommutative power series
$f\in \bbC\langle\langle Z_1, ..., Z_n\rangle\rangle$ lies in
$\Gg_n$ if and only if it is primitive, i.e.,
$\Delta(f) = f\otimes 1 + 1\otimes f$ with respect to the topological
Hopf algebra structure defined above. \

Along with primitive (or Lie) series in $Z_1, ..., Z_n$ we will
consider {\rm group-like elements} of
 $\bbC\langle\langle Z_1, ..., Z_n\rangle\rangle$, i.e., series $\Phi$
 satisfying
 $$\Delta(\Phi) = \Phi\otimes\Phi.\leqno (2.4.6)$$
 The  completed tensor product $\bbC\langle\langle Z_1, ..., Z_n\rangle\rangle
 \widehat{\otimes} \bbC\langle\langle Z_1, ..., Z_n\rangle\rangle$ consists of
 series in $2n$ variables $Z_i'=Z_i\otimes 1$ and $Z_i'' = 1\otimes Z_i$
 which satisfy $[Z_i', Z_j'']=0$ and no other relations. 
 Thus a series $\Phi(Z_1, ..., Z_n)$ is group-like
 if it satisfies the {\em exponential property}: 
$$F(Z_1'+Z_1'', ..., Z_n' + Z_n'') = F(Z_1', ..., Z_n')\cdot F(Z_1'', ..., Z_n''), 
\quad {\rm provided} \quad [Z'_i, Z''_j]=0, \, \forall i,j.\leqno (2.4.7)$$ 
We denote by $G_n$ the set of primitive elements in
$\bbC\langle\langle Z_1, ..., Z_n\rangle\rangle$. Elementary properties of
cocommutative Hopf algebras and elementary convergence  arguments in the
adic topology imply the following:

\proclaim (2.4.8) Proposition. (a) $G_n$ is a group with respect to
the multiplication. \hfill\break
(b) The exponential series defines a bijection
$$\exp: \Gg_n\to G_n,$$
with the inverse given by the logarithmic series.\hfill\break
(c) The image of any series $\Phi\in G_n$ under the commutativization
homomorphism (2.1.8) is a formal series of the form $e^{(a,z)}$
for some $a\in\bbC^n$. \hfill\break
(d) If $\Phi\in G_n$, then 
$$\Phi(-Z_1, ..., -Z_n) = \Phi(Z_1, ..., Z_n)^{-1}$$
(equality of power series).

\vskip .2cm

\noindent {\bf (2.4.9) Example.} 
The above proposition implies that the series 
$$\log(e^{Z_1} \cdot e^{Z_2}) \in {\bbC}\langle\langle Z_1, Z_2\rangle\rangle$$ 
is in fact a Lie series. It is known as the Campbell-Hausdorff series and its 
initial part 
has the form 
$$\log(e^{Z_1}\cdot e^{Z_2}) = Z_1 + Z_2 + {1\over 2} [Z_1, Z_2] + ...$$ 

 \vskip .2cm
 
 Let $G_n(\bbR)\subset G_n$ be the set of group-like series with real coefficients.
 Further, the Lie algebra $FL(Z_1, ..., Z_n)$ is in fact defined over rational
 numbers. In particular, it makes sense to speak about its real part.
 By taking the completion as above, we define the real part of the
 completed free algebra $\Gg_n(\bbR)$.  It is clear
 that the exponential series establishes a bijection between
 $\Gg_n(\bbR)$ and $G_n(\bbR)$. 
 
  The following fact was also pointed out by Chen [C2].
  
  \proclaim (2.4.9) Theorem. If $\gamma\in\Pi_n$ is a path in $\bbR^n$ as above,
  then $E_\gamma(Z)$ is group-like.  Moreover, it lies in the real
  part $G_n(\bbR)$.  
  
  \vskip .1cm
  
  Note that a typical element $\Phi = \Phi(Z_1, ..., Z_n)\in G_n$ is a priori just a formal
  power series and does not have to converge for any given matrix values of the $Z_i$ (unless
  they are all 0). At the same time, series of the form $\Phi=E_\gamma$, $\gamma\in\Pi_n$,
  converge for all values of the $Z_i$. This leads to the proposal, formulated by Chen [C3]
  to view series from $G_n$ with good covergence properties as corresponding to
  ``generalized paths", i.e., paths perhaps more general than piecewise $C^\infty$ ones.
  Theory of
  stochastic integrals, see below, provides a step in a similar direction.

  \vskip .3cm
  
  \noindent {\bf (2.5) Finite-dimensional approximations to $G_n$ and
  $\Gg_n$.} Let us recall a version of the Malcev theory for nilpotent
  Lie algebras. Let $k$ be a field of characteristic 0. A
   Lie algebra $\Gg$ over $k$ is called nilpotent of degree $d$ if all 
$d$-fold iterated commutators in $\Gg$  vanish. Let $U(\Gg)$ be the
universal enveloping algebra of $\Gg$. It is a Hopf algebra with the comultiplication given by
$\Delta(x) = x\otimes 1 + 1\otimes x$ for $x\in \Gg$. 
The subspace $I$ in $U({\Gg})$ generated 
by all nontrivial Lie monomials in elements of $\Gg$, is an ideal, with 
$U({\Gg})/I={\bbC}$. 

\proclaim (2.5.1) Lemma. If $\Gg$ is nilpotent of some degree,
 then $\bigcap I^n = 0$. 
 
 Thus the $I$-adic completion 
$$\widehat{U}({\Gg}) = \lim\limits_{\longleftarrow} U({\Gg})/I^n
\leqno (2.5.2)$$ 
is a complete topological algebra containing $U({\Gg})$. 
As before, the standard Hopf algebra structure on $U({\Gg})$ gives rise to 
a topological Hopf algebra structure on $\widehat{U}({\Gg})$. We then 
have the following fact. 

\proclaim (2.5.3) Theorem. (a) $\Gg$ is the set of primitive elements of 
$\widehat{U}({\Gg})$ \hfill\break 
(b) The set $G$ of group-like elements in $\widehat{U}({\Gg})$ 
is the nilpotent group associated, via the Malcev theory, to
the  Lie algebra $\Gg$. \hfill\break 
(c) If $k = \bbR$ or $\bbC$, then $G$  is the simply connected real or complex
Lie group with Lie algebra $G$. \hfill\break
(d) The exponential map establishes a bijection between $\Gg$ and $G$. 

\vskip .2cm

Let now $k=\bbC$ and
$${\Gg}_{n,d} = FL(X_1, ..., X_n)/FL(X_1, ..., X_n)_{\geq d+1}.
\leqno (2.5.4)$$ 
This is a finite-dimensional Lie algebra known as the free nilpotent Lie algebra 
of degree $d$ generated by $n$ elements. It satisfies the obvious universal 
property. 
 Then 
$${\Gg}_n = \lim_{\longleftarrow}{}_n \,{\Gg}_{n,d}.$$ 
So ${\Gg}_n$ is the free pro-nilpotent Lie algebra on $n$ generators. 

\vskip .2cm

Let $R_{n,d}$ be the quotient of $R_n= {\bbC}\langle\langle Z_1, ..., 
Z_n\rangle\rangle$ 
by the closed ideal generated by all the $(d+1)$-fold commutators of the $Z_i$. 
For example, $R_{n,1} = {\bbC}[[Z_1, ..., Z_n]]$ is the usual  (commutative)
power
series algebra. 

The topological Hopf algebra structure on $R_n$ descends to $R_{n,d}$,
and we easily see the following: 

\proclaim (2.5.5) Proposition.  $R_{n,d}$ is isomorphic to $\widehat{U}({
\Gg}_{n,d})$ 
as a topological Hopf algebra. 

We denote by ${G}_{n,d}\subset R_{n,d}^*$ the group of group-like elements of 
$R_{n,d}$. Then the above facts imply: 

\proclaim (2.5.6) Proposition. (a) ${G}_{n,d}$ is the simply connected complex 
Lie group with Lie algebra ${\Gg}_{n,d}$. \hfill\break 
(b) ${G}_n$ is the projective limit of ${G}_{n,d}$. 

Thus ${G}_{n,d}$ is the ``free unipotent complex algebraic group 
of degree $d$ with $n$ generators" while ${G_n}$ is the free prounipotent 
group. 

As above, taking $k=\bbR$, we get the real parts $G_{n,d}(\bbR)$ and
$\Gg_{n,d}(\bbR)$. The homomorphism $E: \Pi_n\to G_n(\bbR)$
gives rise, for any $d\geq 1$ to the homomorphism
$$\epsilon_{n,d}: \Pi_n\to G_{n,d}(\bbR)\leqno (2.5.7)$$
whose target is a finite-dimensional Lie group. 

\proclaim (2.5.8) Proposition. For any $d\geq 1$ the homomorphism
$\epsilon_{n,d}$ is surjective.

In other words, the group $G_n$ can be seen as a (pro-)algebraic completion
of the path group $\Pi_n$. 

\vskip .2cm

\noindent {\sl Proof:} Let $\Pi_n^{rect}\subset\Pi_n$ be the subgroup of rectangular paths,
i.e., paths consisting of segments each going in the direction of some particular coordinate.
As a group, $\Pi_n^{rect}$ is the free product of $n$ copies of $\bbR$. Let $Z_{i,d}\in\Gg_{n,d}$
be the image of $Z_i$. Then the image of $\Pi_n^{rect}$ in $G_{n,d}(\bbR)$ is the subgroup
generated by the 1-parameter subgroups $\exp (t\cdot Z_{i,d})$, $t\in\bbR$, $i=1, ..., n$.
As the $Z_{i,d}$ generate $\Gg_{n,d}$ as a Lie algebra, the corresponding 1-parameter subgroups
generate $G_{n,d}(\bbR)$ as a group. Therefore $\epsilon_{n,d}(\Pi_n^{rect})=G_{n,d}(\bbR)$.

\vskip .3cm

\noindent {\bf (2.6) Complex exponentials.} Consider the complexification $\bbC^n$ of the space
$\bbR^n$ from (2.1). The form $\Omega$ from (2.1.3) is then a holomorphic form on $\bbC^n$ with
values in $\bbC\langle\langle Z_1, ..., Z_n\rangle\rangle$. In particular, we have the noncommutative
exponential function
$$E_\gamma(Z)\in G_n \subset \bbC\langle\langle Z_1, ..., Z_n\rangle\rangle$$ for any unparametrized
path $\gamma$ in $\bbC^n$ starting at 0. Because $\Omega$ is holomorphic, $E_\gamma(Z)$
is, in addition to invariance under cancellations, also invariant under holomorphic deformations
of sub-paths of $\gamma$. Let $\Pi_n^\bbC$ be the quotient of $\Pi_{2n}$, the group of paths in 
$\bbC^n=\bbR^{2n}$ by the equivalence relation generated by such deformations. Obviously,
$\Pi_n^\bbC$ is a group, and the correspondence $\gamma\mapsto E_\gamma$ gives rise to a 
homomorphism
$$E: \Pi_n^\bbC\to G_n.\leqno (2.6.1)$$
Unlike the real case, it seems to be unknown whether (2.6.1) is injective.
As before, we see that the composite homomorphism
$$\epsilon_{n,d}^\bbC: \Pi_n^\bbC \to G_{n,d}
\leqno (2.6.2)$$
is surjective. 

\vskip .2cm

\noindent {\bf (2.6.3) Example.} Let $C$ be a complex analytic curve, $c_0\in C$ be a point, and
$\phi: C\to \bbC^n$ a holomorphic map such that $\phi(c_0)=0$. Denote by $p: \widetilde{C}\to C$
the universal covering of $C$ corresponding to the base point $c_0$. In other words, $\widetilde{C}$
is ths space of pairs $(c, \gamma)$, where $c\in C$ and $\gamma$ is a homotopy class of paths joinig $c_0$
and $c$. Then, by the above, $\phi$ induces a map $\widetilde{\phi}: \widetilde{C}\to\Pi_n^\bbC$. The
composition
$$\varpi_d = \epsilon_{n,d}^\bbC\circ \widetilde{\phi}: \widetilde{C}\to G_{n,d}$$
can be called the period map of degree $d$. The restriction of $\varpi_d$ to 
$p^{-1}(c_0) = \pi_1(C, c_0)$ is a homomorphism
$$m_d: \pi_1(C, c_0) \to G_{n,d}$$
called the monodromy homomorphism of degree $d$. We get then the ``Albanese map"
$$\alpha_d: C\to G_{n,d}/{\operatorname{Im}}(m_d).$$
The particular case when $C$ is the maximal Abelian covering of a smooth projective curve of
genus $n$, and $\phi$ is the Abel-Jacobi map, corresponds to the setting of Parshin [Pa].
Iterated integrals of modular forms were studied by Manin [Ma]. 

In the subsequent paper [K1] we will use complex noncommutative exponentials to construct
invariants of degenerations of families of curves in an algebraic variety. 

\vfill\eject

\centerline {\bf 3. Generalities on NCFT}

\vskip 1cm

\noindent {\bf (3.0) Formal FT on nilpotent groups.} Let us start with the general situation
of (2.5) with $k=\bbR$. Thus $\Gg$ is a finite dimensional nilpotent real Lie algebra and
$G$ is the corresponding simply connected Lie group. Then $G$ is realized inside $\widehat{U}(\Gg)$
as the set of group-like elements. In general, we can think of elements of $\widehat {U}(\Gg)$ as some
kind of  formal series
(infinite formal linear combinations of elements of a Poincare-Birkhoff-Witt basis of $U(\Gg)$).

 To keep the notation straight, we denote by $E_g\in\widehat{U}(\Gg)$
the element corresponding to $g\in G$. 

\vskip .2cm

\noindent {\bf (3.0.1) Example.} Let $G=\bbR^n$ with coordinates $y_1, ..., y_n$, then $\widehat{U}(\Gg)$
is the ring $\bbC [[z_1, ..., z_n]]$ of formal Taylor series. If $g= (y_1, ..., y_n)\in G$, 
then $E_g = E_g(z) = \exp\left(\sum_i y_i z_i\right)$ is the exponential series with the vector of exponents
$(y_1, ..., y_n)$.

\vskip .2cm

 The above example motivates the following definition. Let $\mu$ be a measure on $G$, or, more generally, a distribution
(understood as a generalized measure, i.e., as a functional on the space of $C^\infty$-functions). Its formal Fourier transform
is the element (formal series) given by
$$\widehat{\cF}(\mu) = \int_{g\in G} E_g d\mu \quad\in\quad \widehat{U}(\Gg), \leqno (3.0.2)$$ 
whenever the integral is defined. 

Recall that for two distributions $\mu,\nu$ on a Lie group $G$ their convolution is defined by
$$\mu*\nu = m_*(\mu\boxtimes\nu),\leqno (3.0.4)$$
where $m: G\times G\to G$ is the multiplication, and $\mu\boxtimes\nu$ is the Cartesian product of
$\mu$ and $\nu$. Here we assume that the pushdown under $m$ is defined. 
The following is then straightforward. 

\proclaim (3.0.5) Proposition. For two (generalized) measures $\mu,\nu$ on $G$ we have
$$\widehat{\cF}(\mu * \nu) = \widehat{\cF}(\mu) \cdot\widehat{\cF}(\nu),$$
(product in $\widehat{U}(\Gg)$). 

\vskip .3cm

\noindent {\bf (3.1) Pro-measures and formal NCFT.}
We now specialize the above to the case when $G = G_{n,d}(\bbR)$. In other words, we
consider the projective system of Lie groups
$$\cdots \to G_{n,3}(\bbR)\to G_{n,2}(\bbR)\to G_{n,1}(\bbR) = \bbR^n\leqno (3.1.1)$$
with projective limit $G_n(\bbR)$. For $d\geq d'$ let
$$p_{dd'}: G_{n,d}(\bbR)\to G_{n,d'}(\bbR)\leqno (3.1.2)$$
be the projection. By a {\em pro-measure} on $G_n(\bbR)$ we will mean
a compatible system of measures on the $G_{n,d}(\bbR)$. In other words, a pro-measure
is a system
${\mu}_\bullet = (\mu_d)$ such that each $\mu_d$ is a measure on
$G_{n,d}(\bbR)$ such that for any $d\geq d'$ the pushdown $(p_{dd'})_*(\mu_d)$
is defined as a measure on $G_{n,d'}(\bbR)$ and is equal to $\mu_{d'}$. Equivalently,
this means that for any continuous function $f$ on $G_{n,d'}(\bbR)$ we have
$$\int_{G_{n,d'}(\bbR)} f\cdot d\mu_{d'} \quad=\quad \int_{G_{n,d}(\bbR)}
(f\circ p_{dd'})\cdot
d\mu_d,\leqno (3.1.3)$$
whenever the LHS is defined. 

\vskip .2cm

More generally, by a pro-distribution we mean a system of distributions on the $G_{n,d}(\bbR)$
(understood as generalized measures, i.e., as functionals on $C^\infty$-functions)
compatible in the similar sense, i.e., satisfying (3.1.3) for $C^\infty$-functions $f$.

\vskip .2cm

For $\Phi = \Phi(Z_1, ..., Z_n)\in G_n$ we denote by $\Phi_{i_1, ..., i_p}$ the coefficient of $\Phi$
at $Z_{i_1} \cdots Z_{i_p}$. It is clear that $\Phi_{i_1, ..., i_p}$ depends only on the image of $\Phi$
in $G_{n,p}$, so it makes sense to speak about $\Psi_{i_1, ..., i_p}$ for $\Phi\in G_{n,d}$, $d\geq p$. 

\vskip .2cm

Let ${\mu}_\bullet$ be a pro-distribution on $G_n(\bbR)$. Its
formal Fourier transform is the formal series $\widehat{\cF}(\mu_\bullet) \in
\bbC\langle\langle Z_1, ..., Z_n\rangle\rangle$ defined as follows:
$$\widehat{\cF}({\mu}_\bullet) = \sum_{p=0}^\infty \sum_{i_1, ..., i_p} \biggl(\int_{\Psi\in G_{n,d}(\bbR)}
\Psi_{i_1 ... i_p} \cdot d\mu_d\biggr) Z_{i_1} \cdots Z_{i_p}.
\leqno (3.1.4)$$
Here for each $p$ the number $d$ is any integer greater or equal to $p$, and we assume that all the integrals converge.

The convolution operation extends, in an obvious way, to pro-distributions on $G_m(\bbR)$ and we get: 

\proclaim (3.1.6) Proposition. If ${\mu}_\bullet, {\nu}_\bullet$ are two
pro-distributions, then
$$\widehat{\cF}({\mu}_\bullet * {\nu}_\bullet) =
\widehat{\cF}({\mu}_\bullet) \cdot \widehat{\cF}({\nu}_\bullet)$$
(product in $\bbC\langle\langle Z_1, ..., Z_n\rangle\rangle$). 

\vskip .3cm

\noindent {\bf (3.2) Delta-functions.} In classical analysis the Fourier transform of $\delta^{(m)}$,
the $m$th derivative of the delta function, is the monomial $z^m$. We now give a noncommutative analog of
this fact.

First of all, let $\delta_d$ be the delta function on $G_{n,d}(\bbR)$ supported at 1. 
Then ${\delta}_\bullet = (\delta_d)$ is a pro-distribution, and
$$\widehat{\cF}({\delta}_\bullet) = 1\in\bbC\langle\langle Z_1, ..., Z_n\rangle\rangle.
\leqno (3.2.1)$$
Next, first derivatives of the delta function at a point on a $C^\infty$ manifold correspond to elements of the
complexified tangent space to the manifold at this point.
This, if $\xi\in FL(Z_1, ..., Z_n)$,
 and $\xi_d$ is the image of $\xi$ in $\Gg_{n,d} = T_1 G_{n,d}(\bbR)\otimes \bbC$, then we have the distribution
 $\partial_{\xi_d}(\delta_d)$ on $G_{n,d}(\bbR)$, and these distributions 
 form a pro-distribution
 $\partial_\xi({\delta}_\bullet)$. 
 
 \vskip .2cm
 
 Further, for any Lie group $G$ with Lie algebra $\Gg$ the iterated derivatives of the delta function at 1
 correspond to elements of $U(\Gg\otimes \bbC)$, the universal enveloping algebra. Thus for any $\psi\in U(\Gg_{n,d})$
we have a punctual distribution $D_\psi(\delta_d)$ on $G_{n,d}(\bbR)$. 

\vskip .2cm

Let now $f\in\bbC\langle Z_1, ..., Z_n\rangle$ be a noncommutative polynomial. Recall that 
$\bbC\langle Z_1, ..., Z_n\rangle$ is the enveloping algebra of $FL(Z_1, ..., Z_d)$. This for any $d$
we have the image of $f$ in $U(\Gg_{n,d})$, which we denote by $f_d$. As before, the distributions $D_{f_d}(\delta_d)$
form a pro-distribution which we denote $D_f({\delta}_\bullet)$. 

\proclaim (3.2.2) Theorem. We have $\widehat{\cF}(D_f({\delta}_\bullet))=f$.
In other words, $\widehat{\cF}$ takes iterated derivatives of the delta function into
(noncommutative) polynomials.

\vskip .2cm

Let $L_{f,d}$ be the left invariant differential operator on $G_{n,d}(\bbR)$ corresponding to $f_d\in U(\Gg_{n,d}))$.
Similarly, let $R_{f,d}$ be the right invariant differential operator corresponding to $f_d$. Recall that
distributions (volume forms) form a right module over the ring of differential operators. In other words, if $P$
is a differential operator acting on functions by $\phi\mapsto P\phi$, then we write the action of the
adjoint operator on volume forms by $\omega\mapsto \omega P$. 
Thus,  if
${\mu}_\bullet = (\mu_d)$ is a pro-distribution, and $f\in \bbC\langle Z_1, ..., Z_n\rangle$,
then we have pro-distributions ${\mu}_\bullet L_f = (\mu_d L_{f,d})$ and
${\mu}_\bullet R_f = (\mu_d R_{f,d})$. Since applying $R_{f,d}$ or $L_{f,d}$ to a distribution
is the same as the right or left convolution with $D_{f_d}(\delta_d)$, Proposition 3.1.6 implies the following.

\proclaim (3.2.3) Proposition. If $\phi\in \bbC\langle\langle Z_1, ..., Z_n\rangle\rangle$ is the
Fourier transform of ${\mu}_\bullet$, then for any $f\in\bbC\langle Z_1, ..., Z_n\rangle$
the product $f\cdot \phi$ is the Fourier transform of ${\mu}_\bullet L_f$, and
$\phi\cdot f$ is the Fourier transform of ${\mu}_\bullet R_f$.

\vskip .3cm

\noindent {\bf (3.3) Measures and convergent NCFT.} Let $p_d: G_n(\bbR)\to G_{n,d}(\bbR)$ be the
projection. By a cylindrical open set in $G_n(\bbR)$ we mean a set of the form $p_d^{-1}(U)$,
where $d\geq 1$ and $U\subset G_{n,d}(\bbR)$ is an open set. These sets form thus a basis of the
projective limit topology on $G_n(\bbR)$. We denote by ${\frak S}$ the $\sigma$-algebra of sets in
$G_n(\bbR)$ generated by cylindrical open sets. Its elements will be simply called Borel subsets
in $G_n(\bbR)$. 

\vskip .2cm

\noindent {\bf (3.3.1) Example.} Let
$$G_n(\bbR)^{ent} = G_n(\bbR)\cap \bbC\langle\langle Z_1, ..., Z_n\rangle\rangle ^{ent}$$
be the subgroup formed by entire series, see (1.2.2). Since for $\Phi\in G_n(\bbR)$
each given coefficient of $f$ depends on the image of $\Phi$ in some $G_{n,d}(\bbR)$, the condition
(1.2.2) implies that $G_n(\bbR)^{ent}$ is a Borel subset. Note further that for 
$\Phi\in G_n(\bbR)^{ent}$ and any Hermitian matrices $Z_1, ..., Z_n$ (of any size $N$)
the matrix $\Phi (iZ_1, ..., iZ_n)$ is unitary. This follows from the reality of the
coefficients in $\Phi$ and from Proposition 2.4.8(d). 

\vskip .2cm

By a measure on $G_n(\bbR)$ we mean a complex valued, countably additive measure on the
$\sigma$-algebra ${\frak S}$. If $\mu$ is such a measure, we define its Fourier transform
to be the function of indeterminate Hermitian $N$ by $N$ matrices $Z_1, ..., Z_n$
(with indeterminate $N$) given by
$$\cF(\mu) (Z_1, ..., Z_n) = \int_{\Phi\in G_n(\bbR)^{ent}} \Phi(iZ_1, ..., iZ_n) d\mu(\Phi).\leqno
(3.3.2)$$
As usual, by a probability measure on $G_n(\bbR)$ we mean a real, nonnegative-valued measure
on ${\frak S}$ of total volume 1.  
\vskip .2cm

Given a pro-measure ${\mu}_\bullet = (\mu_d)$ on $G_n(\bbR)$, the correspondence
$$p_d^{-1}(U)\mapsto \mu_d(U), \quad U\in G_{n,d}(\bbR), 
\leqno (3.3.3)$$
defines a finite-additive function on cylindrical open sets in $G_n(\bbR)$. The following
fact is a version of the basic theorem of Kolmogoroff ( [SW], Thm. 1.1.10) that a stochastic process is uniquely determined
by its finite-dimensional distributions. More precisely, Kolmogoroff's theorem is about  probability measures on
an infinite product of measure spaces. The modification to the case of  a projective limit is immediate.

\proclaim (3.3.4) Theorem. If ${\mu}_\bullet$ is a probability pro-measure (i.e., each
$\mu_d$ is a probability measure), then the correspondence (3.3.4) extends to a unique probability measure
$\mu = \lim\limits_{\longleftarrow} \, \mu_d$ on $G_n(\bbR)$, so that $\mu_d = p_{d*}(\mu)$.

Thus, probability measures and pro-measures are in bijection. 

\vfill\eject

\centerline{\bf 4. Noncommutative Gaussian and the Wiener measure}

\vskip 1cm

\noindent {\bf (4.1) Informal overview.} By the noncommutative Gaussian we mean the following noncommutative
power series
$$\Xi(Z) = \exp\biggl( -{1\over 2} \sum_{i=1}^n Z_i^2\biggr) \quad\in\quad \bbC
\langle\langle Z_1, ..., Z_n\rangle\rangle^{ent}.\leqno (4.1.1)$$
As the series is entire, we will denote by the same symbol $\Xi(Z_1, ..., Z_n)$
its value on any given square matrices $Z_1, ..., Z_n$. In the classical (commutative) analysis
the Fourier transform of a Gaussian is another Gaussian. In this section we present a noncommutative
extension of this fact. Informally the answer can be formulated as follows.

\proclaim (4.1.2) Informal theorem. ``The" measure on the space of paths whose Fourier transform
gives $\Xi(Z)$, is the Wiener measure.

\vskip .2cm

We write ``the" in quotes because so far there is no uniqueness result for NCFT, so (4.1.2) 
can be read in one direction: that the NCFT of the Wiener measure is $\Xi(Z)$. Still, 
there are two more issues one has to address in order to make (4.1.2) into a theorem. 
First, the
Wiener measure (see below for a summary) is defined on the space of parametrized paths, while
NCFT is defined for measures on the space of unparametrized paths. This can be addressed by considering
the pushdown of the Wiener measure (i.e., by performing the
integration over the space of parametrized paths).

Second, and more importantly, the Wiener measure is defined on the space of continuous paths,
and piecewise smooth paths form a subset of measure 0. On the other hand, the series $E_\gamma(Z)$
is a solution of a differential equation involving the time derivatives of $\gamma$ and so is
a priori not defined if $\gamma$ is just a continuous path. This difficulty is resolved by using the
theory of stochastic integrals and stochastic differential equations which indeed provides
a way of associating $E_\gamma(Z)$ to all continuous $\gamma$ except those forming a set of Wiener measure
0. 

Once these two modifications are implemented, (4.1.2) becomes an instance of  the familiar
principle in the theory of stochastic differential equations: that the direct image of
the Wiener measure under the map given by the solution of a stochastic differential equation,
is the heat measure for the corresponding (hypo)elliptic operator, see [Bel] [Ok] [Bi1].

\vskip .3cm

\noindent {\bf (4.2) The hypo-Laplacians and their heat kernels.} Let $Z_{i,d}$
be the image of $Z_i$ in $\Gg_{n,d}$, and $L_{i,d}$ be the
left invariant vector field on $G_{n,d}(\bbR)$ corresponding to $Z_{i,d}$. We consider
$L_{i,d}$ as a first order differential operator on functions. The $d$th hypo-Laplacian is the
operator
$$\Delta_d  = \sum_{i=1}^n L_{i,d}^2\leqno (4.2.1)$$
in functions on $G_{n,d}(\bbR)$. For $d\geq d'$ the operators $\Delta_d$ and $\Delta_{d'}$
are compatible:
$$\Delta_d(p_{dd'}^*f) = p_{dd'}^* (\Delta_{d'}f), \quad\forall f\in C^2(G_{d'}(\bbR)),
\leqno (4.2.2)$$
where the projection $p_{dd'}$ is as in (3.1.2). This follows because a similar compatibility
holds for each $L_{i,d}$ and $L_{i, d'}$. 
 
For $d>1$ the number of summands in (4.2.1) is less than the dimension of $G_{n,d}(\bbR)$,
so $\Delta_d$ is not elliptic. However, $\Delta_d$ is hypoelliptic [Ho], i.e.,
every distribution solution of $\Delta_d u=0$ is real analytic.
This follows from Theorem 1.1 of H\"ormander [Ho], since the $Z_{i,d}$ generate $\Gg_{n,d}$
as a Lie algebra. Further, it is obvious that $\Delta_d$ is positive:
$$(\Delta_d u,u)\geq 0, \quad u\in C^\infty_0(G_{n,d}(\bbR)).
\leqno (4.2.3)$$
General properties of positive hypoelliptic operators [Ho] imply that the heat
operator $\exp(-t \Delta_d)$, $t>0$, is given by a positive $C^\infty$ kernel. Because this
operator is left invariant, we get part (a) of the following theorem:

\proclaim (4.2.4) Theorem. (a) The operator $\exp(\Delta_d/2)$ is given by convolution with a uniquely
defined probability measure $\theta_d$ on $G_{n,d}(\bbR)$. This measure is real analytic with
respect to the Haar measure. \hfill\break
(b) For $d\geq d'$ the measures $\theta_d$ and $\theta_{d'}$ are compatible:
$(p_{dd'})_*(\theta_d) = \theta_{d'}$.

Part (b) above follows from (4.2.2). 

\vskip .2cm

Thus we obtain a probability pro-measure ${\theta}_\bullet = (\theta_d)$
on $G_n(\bbR)$ and hence a probability measure $\theta = \lim\limits_{\longleftarrow} \theta_d$. 

\vskip .2cm

\noindent {\bf (4.2.5) Examples.} (a) the group $G_{n,1}$ is identified with the space $\bbR^n$
from (2.1) with coordinates $y_1, ..., y_n$, and $Z_{i,1}= \partial/\partial y_i$. Therefore
$\Delta_1$ is the standard Laplacian on $\bbR^n$, and 
$$\theta_1 = {dy_1 ... dy_n\over (2\pi)^{n/2}} \exp\biggl(-{1\over 2} \sum_{i=1}^n y_i^2\biggr)$$
is the usual Gaussian measure on $\bbR^n$. Each $\theta_d$, $d>1$ is thus a lift of this measure
to $G_{n,d}$.

\vskip .1cm

(b) For $d=2$ an explicit formula for $\theta_2$ was obtained by Gaveau in [G].
Here we consider the case $n=2$  where the formula was also obtained by Hulanicki [Hu].
 In this case $\Gg_{2,2}$ is the Heisenberg Lie algebra 
with basis consisting of $Z_{1,2}, Z_{2,2}$ and the 
central element $h=[Z_{1,2}, Z_{2,2}]$. Denoting $y_1, y_2, v$ the corresponding
exponential coordinates on $G_{2,2}$, we have
$$\theta_{2} = {dy_1dy_2 dv \over (2\pi)^2} \int_{\tau=-\infty}^\infty {2\tau \over \sinh (2\tau)} \cdot
\exp\biggl( i\tau v - (y_1^2+y_2^2) {2\tau\over \tanh (2\tau)} \biggr)d\tau.$$
In fact, all known formulas in the literature (see [BGG] for a survey) involve integration over auxuliary parameters.

\vskip .2cm

\proclaim (4.2.6) Theorem. The formal Fourier transform of
the pro-measure ${\theta}_\bullet$ is equal to the noncommutative Gaussian
$\Xi(Z)$. 

\noindent {\sl Proof:} This follows from the fact that the delta-pro-distribution $D_{Z_i}(\delta_\bullet)$ corresponding to
the generator $Z_i\in \Gg_n$ is taken by $\cF$ into the monomial $Z_i$. For each $d$ the corresponding
distribution takes a function $f$ on $G_{n,d}(\bbR)$ into the value of $L_i(f)$ at the unit element of
$G_{n,d}(\bbR)$. Further, 
convolution of such
pro-distributions corresponds to composition of left invariant differential operators in
the spaces of functions of the $G_{n,d}(\bbR)$. So the system of the heat kernel operators
on the $G_{n,d}(\bbR)$, $d\geq 1$, given by $\exp\left(-{1\over 2}\sum L_{i}^2\right)$
has, as a pro-distribution, the Fourier transform equal to $\exp\left(-{1\over 2}\sum Z_i^2\right)$. \qed

\vskip .3cm

\noindent {\bf (4.3) The Wiener measure.} Let $P_n$ be the space of
continuous parametrized paths $\gamma: [0,1]\to \bbR^n$ such that $\gamma(0)=0$.
The Wiener measure $w$ on $P_n$ is first defined on cylindrical open sets
$C(t_1, ..., t_m, U_1, ..., U_m)$, where $0<t_1 < ..., <t_m < 1$ and $U_i\subset \bbR^n$ is open.
By definition,
$$C(t_1, ..., t_m, U_1, ..., U_m) = 
\bigl\{ \gamma: \,\, \gamma(t_i)\in U_i, \,\,i=1, ..., m\bigl\},$$
and
$$w\bigl(C(t_1, ..., t_m, U_1, ..., U_m)\bigr) = \leqno (4.3.1) $$
$$=\int_{(y^{(1)}, ..., y^{(m)})\in U_1\times ...\times U_m}
\prod_{i=0}^m {\exp\bigl(-\|y^{(i+1)}-y^{(i)}\|^2/2(t_{i+1}-t_{i})\bigr)
\over \bigl(2\pi (t_{i+1}-t_i)\bigr)^{1/2}}
dy^{(1)} ... dy^{(m)}.$$
Here we put $t_0=0, t_{m+1}=1$ and $y^{(0)}=0$. Further, it is proved that
$w$ extends to a probability
measure on the $\sigma$-algebra generated by the above subsets.

\vskip .2cm

The Brownian motion is the family of $\bbR^n$-valued functions
(random variables) on $P_n$ parametrized by $t\in [0,1]$:
$$b(t) = (b_1(t), ..., b_n(t)), \quad b(t): P_n\to \bbR^n, \,\, b(t)(\gamma) =
\gamma(t). \leqno (4.3.2) $$
let $P_n^{sm}\subset P_n$ be the subset of piecewuse smooth paths. Then it is well known that $w(P_n^{sm})=0$. 

As well known, the Wiener measure has the following intuitive interpretation
$$dw(\gamma)  = \exp\left( -\int_0^1 \|\gamma'(t)\|^2 dt\right) \cD\gamma, \quad \cD\gamma = \prod_{t=0}^1 d\gamma(t). \leqno (4.3.3)$$
In other words, $\cD\gamma$ is the (nonexistent) Lebesgue measure on the infinite-dimensional vector space of all paths,
while the integral in the exponential is the action of a free particle.

\vskip .3cm

\noindent {\bf (4.5) Reminder on stochastic integrals.} Let $\omega = \sum_{i=1}^n \phi_i(y) dy_i$ be a
1-form on $\bbR^n$ with (complex valued) $C^\infty$ coefficients. If $\gamma: [0,1]\to \bbR^n$
is a piecewise smooth path, then we can integrate $\omega$ along $\gamma$, getting a number
$$\int_\gamma \omega = \int_0^1 \gamma^*(\omega) = \int_0^1 \sum_i \phi_i(\gamma(t))\gamma_i'(t) dt.
\leqno (4.5.1)$$
This gives a map (function)\
$$\int (\omega): P_n^{sm}\longrightarrow \bbC. \leqno (4.5.2) $$
If $\gamma(t)$ is just a continuous path without any differentiability
assumptions, then (4.5.1) is not defined, so there is no immediate extension of the map (4.5.2)
to the space $P_n$. The theory of stochastic integrals provides several (a priori different) ways 
to construct such an extension. The two best known approaches are the Ito and Stratonovich integrals over the
Brownian motion, see [SW] [KW]. They are functions
$$\int^{\Ito} (\omega), \quad \int^{\Str} (\omega): \quad P_n\to \bbC, \leqno (4.5.3)$$
defined everywhere outside some subset of Wiener measure 0, and measurable with respect to this measure. 
To construct them, see [Ok], p. 14-16, one has to consider Riemann sum approximations to the integral but
restrict to Riemann sums of some particular type. For a piecewise smooth path $\gamma$
the integral is the limit of sums
$$\sum_{i=1}^n \sum_{\nu=1}^m \phi_i(\gamma(\xi_\nu))\bigl( (\gamma_i(t_\nu) - \gamma_i(t_{\nu -1})\bigr), 
\leqno (4.5.4)$$
where $0 = t_0 < t_1 < ... < t_m=1$ is a decomposition of $[0,1]$ into intervals, and $\xi_\nu\in [t_{\nu -1}, t_\nu]$
are some chosen points. In the smooth case the limit exists provided $\max (t_\nu - t_{\nu -1})$ goes to 0
(in particular, the choice of $\xi_\nu$ is inessential). 

Now, to obtain $\int^{\Ito} (\omega)$, one chooses the class of Riemann sums with
$$t_\nu = \nu/m, \quad \xi_\nu = t_{\nu -1}, \quad m= 2^q, q\to\infty. \leqno (4.5.5)$$
 In other words, for each $q$ the above sum defines a function
$\Sc_q^{\Ito} (\omega): P_n\to \bbC$, and
$$\int^{\Ito} (\omega) = \lim_{q\to\infty} \Sc_q^{\Ito} (\omega). \leqno (4.5.6)$$
To obtain $\int^{\Str}(\omega)$, one chooses the class of Riemann sums with
$$t_\nu = \nu/m, \quad \xi_\nu = (t_{\nu -1} + t_\nu)/2, \quad m=2^q, q\to\infty. \leqno (4.5.7)$$
Each such sum gives a function $\Sc_q^{\Str} (\omega): P_n\to \bbC$, and
$$\int^{\Str}(\omega) = \lim_{q\to\infty} \Sc_q^{\Str}(\omega).\leqno (4.5.8)$$ 
It follows that $\int^{\Str}(\omega)$ is invariant under smooth reparametrizations of the path
considered as transformations acting on $P_n$.

The more common notation for the stochastic integrals (considered as random variables on $P_n$) is:
$$\int^{\Ito}(\omega) = \int_0^1 \omega(b(t)) db(t)), \quad \int^{\Str}(\omega) = \int_0^1 \omega(b(t)) \, \circ db(t)), 
\leqno (4.5.9)$$
where $b(t)$ is the Brownian motion (4.3.2). Thus $db(t)$ and $\circ db(t)$ stand for the two ways
(due to Ito and Stratonovich) of regularizing the (a priori divergent) differential of the Brownian path $b(t)$. 
See [Ok] for the relation between the two regularization schemes. 
By restricting to the truncated path $[0,s]$, $s\leq 1$, one defines the stochastic integrals $\int_0^s$ in each of the
above setting. 

\vskip .3cm

\noindent {\bf (4.6) Stochastic holonomy.} Let $G$ be a Lie group which we suppose to be embedded as
a closed subgroup of $GL_N(\bbC)$ for some $N$ and let $\Gg \subset \Mat_N(\bbC)$ be the Lie algebra of $G$.
Let $A=\sum A_i(y) dy_i$ be a smooth $\Gg$-valued 1-form on $\bbR^n$, which we consider as a connection
in the trivial $G$-bundle over $\bbR^n$. If $\gamma: [0,1]\to \bbR^n$ is a piecewise smooth
path, then we have the holonomy of $A$ along $\gamma$: 
$$\Hol_\gamma(A) = P\exp \int_\gamma A \quad \in\quad G.\leqno (4.6.1)$$
It is the value at $t=1$ of the solution $U(t)\in GL_N(\bbC)$ of the
differential equation
$${dU\over dt} = U(t) \biggl( \sum_i A_i(\gamma(t)) \cdot \gamma_i'(t)\biggr), \quad U(0)=1. \leqno (4.6.2)$$
The holonomy defines thus a map
$$\Hol (A): P_n^{sm}\to G.\leqno (4.6.3)$$
As before, (4.6.2) and thus $\Hol_\gamma(A)$ have no immediate sense without some differentiability
assumptions on $A$. 
The theory of stochastic differential equations [Ok] [KW] resolves this difficulty by replacing the above
differential equation by an integral equation and understanding the integral in a regularized sense as in (4.5).
Thus, one defines the  Ito and Stratonovich  stochastic holonomies  which are measurable maps
$$\Hol^{\Ito}(A), \,\,\Hol^{\Str}(A): P_n\to G, \leqno (4.6.3)$$
defined outside a subset of Wiener measure 0. For example, $\Hol^{\Str}(A)$ is defined as the value at $t=1$
of the $G$-valued stochastic process $U(t)$ satisfying the Stratonovich integral equation
$$U(t) = 1 + \int_0^t  U(s)  \biggl( \sum_i A_i(b(s)) \circ db_i(s)\biggr). \leqno (4.6.4)$$
We will be particularly interested in the case when $A_i$ are constant, i.e., our connection is translation
invariant. In this case $B(t) = \sum A_i b_i(t)$ is a (possibly degenerate) Brownian motion on $\Gg$
and $U(t)$ is the corresponding left invariant Brownian motion on $G$ as studied by McKean,  see [McK], \S 4.7,
and also [HL]. In particular, the Stratonovich holonomy can be represented as a ``product integral"
in the sense of McKean
$$\Hol^{\Str}(A) = \prod_{t\in [0,1]} \exp(dB(t)) \quad := 
\quad \lim_{q\to\infty} \,\,\prod_{\nu=1}^{2^q} 
\exp\biggl( B\left({\nu\over 2^q}\right) - B\left({\nu - 1 \over 2^q}\right)\biggr), \leqno (4.6.5)$$ 
see [HL], Thm.2. Here the product is taken in the order of increasing $\nu$. 

Again, the product integral representation implies invariance of $\Hol^{\Str}(A)$ under smooth
reparametrizations of the path. 

\vskip .3cm

\noindent {\bf (4.7) The  Malliavin calculus and the Feynman-Kac-Bismut formula.} 
We now specialize (4.6) to the case when $G= G_{n,d}(\bbR)$, $\Gg = \Gg_{n,d}(\bbR)$ and
$A= A^{(d)}$ is the constant 1-form $A^{(d)} = \sum_{i=1}^n Z_{i,d} dy_i$. 
We get the stochastic holonomy map
$$\Hol^{\Str}(A^{(d)}): P_n\to G_{n,d}(\bbR). \leqno (4.7.1)$$

\proclaim (4.7.2) Theorem. The probability measure $\theta_d$ on $G_{n,d}(\bbR)$ is equal to
$\Hol^{\Str}(A^{(d)})_*(w)$, the push-down of the Wiener measure under the holonomy map. 

\noindent {\sl Proof:} This is a fundamental property of (hypo)elliptic diffusions holding for any
vector fields $\xi_1, ..., \xi_n$ on a manifold $M$ such that iterated commutators of the $\xi_i$
span the tangent space at every point. In this case the operator $\Delta = \sum \Lie_{\xi_i}^2$
is hypoelliptic and  has a uniquely defined, smooth heat kernel $\Theta (x,y), x,y\in M$ 
which is a function in $x$
and a volume form in $y$ and represents the operator $\exp(-\Delta/2)$. Further, the heat equation
$$\partial u/\partial t = -\Delta(u)/2\leqno (4.7.3)$$ is the ``Kolmogoroff backward equation" for the $M$-valued stochastic 
process $U(t)$ satisfying the Stratonovich differential equation
$$dU = \sum L_{\xi_i}(U) \circ db_i\leqno (4.7.4)$$
with the $b_i(t)$ being as before. This means that the fundamental solution of (4.7.3) is
the pushforward of the Wiener measure under the process $U(t)$. See [Ok], Th. 8.1.
Our case is obtained by specializing to $M=G_{n,d}(\bbR)$, $\xi_i = Z_{i,d}$. \qed

\vskip .2cm

Further, let $\theta$ be the probability measure on 
$G_n(\bbR) = \lim\limits_{\longleftarrow}{}_d \,G_{n,d}(\bbR)$
corresponding to the pro-measure $(\theta_d)$ by Theorem 3.3.5.
 Note that the maps $\Hol^{\Str}(A^{(d)})$ for various $d$
unite into a map
$$\Hol^{\Str}(A): P_n\to G_n(\bbR), \quad A = \sum Z_i dy_i.  
\leqno (4.7.4)$$
We get the following corollary.

\proclaim (4.7.5) Corollary. The measure $\theta$ is the pushdown of the Wiener measure under $\Hol^{\Str}(A)$.

\vskip .2cm

\proclaim (4.7.6) Theorem. (a) The support of the measure $\theta$ is 
contained in $G_n(\bbR)^{ent}$,
the set of entire group-like power series. \hfill\break
(b) The convergent Fourier transform of $\theta$ is equal to $\Xi(Z)$. In other words
(taking into account part (a) and (4.7.5)),  for any
given Hermitian matrices $Z_1, ..., Z_n$ of any given size $N$ we have
$$\exp\biggl(-{1\over 2} \sum_{j=1}^n Z_j^2\biggr) = \int_{\gamma\in P_n} \Hol^{\Str}_\gamma (A) (iZ_1, ..., iZ_n) dw(\gamma).$$

\noindent {\sl Proof:} (a) First of all, one can write $\Hol^{\Str}(A)$,  similarly to (2.3.2), as the generating function
of stochastic  iterated integrals understood in the sense of Stratonovich:
$$\Hol^{\Str}(A) = \sum_m \sum_{J= (j_1, ..., i_m)} Z_{j_1} ... Z_{j_m} \cdot \int^{\rightarrow} \circ db_{j_1} ... \circ db_{j_m}.
\leqno (4.7.7)$$
See [FN] for the definition of stochastic iterated integrals as well as for the proof. So the question is whether the
series (4.7.7)  (with coefficients being random variables on $P_n$) is entire almost surely (i.e., outside a set of
Wiener measure 0). Questions of this nature (``convergence of stochastic Taylor series") were studied by Ben Arous [Be]. 
To get a link to his notation, 
denote the stochastic integral in the RHS of (4.7.7) by $B_J$ and the monomial $Z_{j_1}... Z_{j_m}$ by
simply $Z^J$. The stochastic Taylor series of [Be] are expressions of the
form
$$\sum_m \sum_{j=(j_1, ..., j_m)} x_J \cdot B_J, \leqno (4.7.8)$$
where $(x_J)$ is a family of (nonrandom) complex numbers labelled by the multi-indices $J$. 
Let us write $m=l(J)$ for the total degree of the monomial corresponding to $J$. 
For the noncommutative formal series $\Hol^{\Str}(A) = \sum_J B_J Z^J$ to be entire (1.2.2), it should satisfy
$|B_J| \cdot R^{l(J)}\to 0$ for any $R>0$ which is equivalent to saying that $\sum_J |B_J|\cdot R^{l(J)} <\infty$
for any $R>0$. The last sum is an example of a Ben Arous series, and Corollary 1 of [Be] gives its
almost sure convergence, whence the claim.

\vskip .1cm

(b) It follows from (a) and Theorem 4.2.6 about the formal Fourier transform. \qed

\vfill\eject

\centerline {\bf 5. Futher examples of NCFT}

\vskip 1cm

\noindent {\bf (5.1) Near-Gaussians.} In classical analysis, a near-Gaussian is a function of the form
$f(z)\cdot e^{-{\|z\|}^2 /2}$ where $f(z)$, $z=(z_1, ..., z_n)$, is a polynomial. In that setting, the Fourier
transform of a near-Gaussian is another near-Gaussian. 

A natural noncommutative analog of a near-Gaussian is a function of the form
$$ F(Z) \cdot \Xi(Z) \cdot G(Z), \quad F, G \in \bbC\langle Z_1, ..., Z_n\rangle. \leqno (5.1.1)$$
 It can be represented as a (formal) Fourier transform using Proposition 3.2.3:
 $$ F(Z)\cdot \Xi(Z)\cdot G(Z) =  \widehat{\cF}\bigl( \theta_\bullet L_F R_G), \leqno (5.1.2)$$
 where $L_F$ is the system of left invariant differential operators on the $G_{n,d}(\bbR), d\geq 1$,
 corresponding to $F$, while $R_G$ is the system of right invariant differential operators corresponding to
 $G$. 
 
 It seems difficult to realize the measures $\theta_d L_F R_G, \quad d\geq 1$, in terms of some transparent measures
 on the space $P_n$, as it requires using group translations on $\Pi_n^{cont}$, the group of
 continuous paths obtained by quotienting $P_n$ by reparametrizations and cancellations.

 \vskip .3cm
 
 \noindent {\bf (5.2) The Green pro-measure.} Let $g_d$ be the fundamental solution of
 the $d$th hypo-Laplacian on $G_{n,d}(\bbR)$ centered at $1$, the unit element, i.e.,
 $$\Delta_d (g_d) = \delta_1. \leqno (5.2.1)$$
 By the general properties of hypoelliptic operators, $g_d$ is a measure (volume form)
 on $G_{n,d}(\bbR)$ smooth away from $1$. In fact, if we denote by $\theta_{d,t}$ the kernel
 of $\exp(-t\Delta_d/2)$, $t>0$, i.e., the heat kernel measure at time $t$, then
 $$g_d = \int_{t=0}^\infty \theta_{d,t} dt. \leqno (5.2.2)$$
 This expresses the fact that the Green measure of a domain is equal to the
 amount of time a diffusion path spends in the domain. 
 It is clear therefore that $g_\bullet = (g_d)$ is a pro-measure on $G_n(\bbR)$.
 
 \vskip .2cm
 
 \noindent {\bf (5.2.3) Examples. } (a) For $d=1$ we have the Green function of the usual
 Euclidean Laplacian in $\bbR^n$ which has the form
 $$g_1(y) = {1\over 4\pi} \ln (y_1^2 + y_2^2)dy_1 dy_2, \quad n=2,$$
 $$g_1(y) = -{((n/2)-2)!\over 4 \pi^{n/2}} \left(\sum y_i^2\right)^{1-n/2} dy_1 ... dy_n. $$ 
 
 (b) Consider the case $n=2, d=2$ corresponding to the Heisenberg group, and let us use the
 exponential coordinates $y_1, y_2, v$ as in Example 4.2.5(b). Then
 $$g_2 (y_1, y_2, v) = {1\over\pi} {1\over \sqrt{(y_1^2+ y_2^2) + v^2}} \, dy_1dy_2 dv,$$
 as was found by Folland [Fo],  see also [G], p. 101. 
 
 \vskip .3cm
 
 \noindent {\bf (5.3) The method of kernels.}  More generally, if $F(Z_1, ..., Z_n)$
 is a ``noncommutative function"
 such that the operator $F(L_{1,d}, ..., L_{n,d})$ in functions on $G_{n,d}(\bbR)$ makes sense and posesses a 
 distribution kernel
 $K_d(x,y)dy$, then the distribution $\mu_d = K_d(1,y)dy$ is precisely the $d$th component
 of the pro-distribution whose Fourier transform is $F$. 
 
 For example, hypoelliptic calculus allows us to consider $F(Z) = \phi\left(\sum_{i=1}^n Z_i^2\right)$
 where $\phi: \bbR\to \bbR$ is any $C^\infty$ function decaying at infinity such as $\phi(u) = e^{-u^2/2}$
 or $\phi(u) = 1/u$, or $1/(u^2+1)$. This leads to a considerable supply of pro-distributions.

\vskip .3cm

\noindent {\bf (5.4) Probabilistic meaning.} An idea in probability theory very similar to our NCFT, viz. the idea
of associating a
noncommutative power series to a stochastic process, was proposed by Baudoin [Ba], who called this
series ``expectation of the signature" and emphasized its importance. From the general viewpoint of
probability theory one can look at this series (the Fourier transform of a probability measure on the
space of paths) as being rather an analog of the characteristic function of $n$ random variables.
Indeed, if $x_1, ..., x_n$ are random variables, then their joint distribution is a probability measure
on $\bbR^n$, and the characteristic function is the (usual) Fourier transform of this measure:
$$f(z_1, ..., z_n) = {\Bbb E}\bigl[e^{i (z,x)}\bigr] \leqno (5.4.1)$$
which is an entire function of $n$ variables. 
Each time when we have a natural lifting of the characteristic function to the noncommutative domain,
we can therefore expect some $n$-dimensional stochastic process lurking in the background.  

\vfill\eject

\centerline {\bf 6. Fourier transform of noncommutative measures}

\vskip 1cm

\noindent {\bf (6.1) Nomcommutative measures.} Following the general approach of Noncommutative Geometry [Con],
we consider a possibly noncommutative $\bbC$-algebra $R$ (with unit) as a replacement of a ``space"
(${\operatorname{Spec}}(A)$). A measure on $R$ is then simply a linear functional (``integration map")
$\tau: I\to\bbC$ defined on an appropriate subspace $I\subset R$ whose elements have the meaning 
of integrable functions. We will call a measure $\tau$ {\em finite}, if $I=R$,
and
{\em normalized}, if it is finite and $\tau(1)=1$. If $R$ has a structure of a $*$-algebra, then
a finite measure $\tau$ is called positive if $\tau(aa^*)\geq 0$ for any $a\in R$.
A (noncommutative) {\it probability measure} on a $*$-algebra $A$ is a normalized, positive measure.

\vskip .2cm

\noindent {\bf (6.1.1) Examples.} (a) Let $R= \Mat_N(\bbC)$ with the $*$-algebra structure given by
the Hermitian conjugation. Then $\tau(a) = {1\over N} \Tr(a)$ is a probability measure. 

\vskip .1cm

(b)  Let $R = \bbC\langle Z_1, ..., Z_n\rangle$ with the  $*$-algebra structure given by $Z_i^*=Z_i$.
Let $\Herm_N$ be the space of Hermitian $N$ by $N$ matrices. We denote by $dZ= \prod_{i,j=1}^N  dZ_{ij}$
the standard volume form on $\Herm_N$. 
Let $\mu= \mu_N$ be a volume form on $(\Herm_N)^n$ of exponential decay at infinity. Then we have a finite 
 measure on $R$
given by
$$\tau(f) = {1\over N} \Tr \int_{Z_1, ..., Z_n\in\Herm_N} f(Z_1, ..., Z_n) d\mu(Z_1, ..., Z_n). $$
If $\mu_N$ is a normalized (resp. probability) measure in the usual sense, then $\tau$ is a 
normalized (resp. probability) measure in the
noncommutative sense. An important example is
$$\mu_N = \exp\left( -S(Z_1, ..., Z_n)\right) dZ_1 ... dZ_n,$$
where the ``action"  $S(Z_1, ..., Z_n)\in \bbC\langle Z_1, ..., Z_n\rangle$ is a noncommutative
polynomial with appropriate growth conditions at infinity. 

\vskip .1cm

(c) Let $R = \bbC\langle X_1^{\pm 1}, ..., X_n^{\pm 1}\rangle$ with the $*$-algebra structure given by
$X_i^* = X_i^{-1}$. If $\mu = \mu_N$ is a finite measure on $U(N)^n$, then we have a  finite measure $\tau$ on $R$
given by
$$\tau(f) = {1\over N} \Tr \int_{X_1, ..., X_n\in U(N)} f(X_1, ..., X_n) d\mu(X_1, ..., X_n),$$
which is normalized (resp. probability) if $\mu_N$ is so in the usual sense.

\vskip .3cm

\noindent {\bf (6.2) Free products.} Let $R_1, ..., R_n$ be algebras with unit. Then we have their
free product $R_1\star .. \star R_n$.  
This is an 
algebra 
containing all the $R_i$ and characterized by the following universal property: 
for any algebra $B$ and any homomorphisms $f_i: R_i\to B$ there is a unique 
homomorphism $f: R_1\star ...\star R_n\to B$ restricting to $f_i$ on $R_i$ for each $i$. 
Explicitly, $R_1\star ...\star R_n$ is obtained as the quotient of the free (tensor) 
algebra generated by the vector space $R_1\oplus ... \oplus R_n$ by the 
relations saying the  products of elements from each $R_i$ are given 
by the existing multiplication in $R_i$. We will also use the notation
$\bigstar_{i=1}^n R_i$.

\vskip .1cm 

\noindent {\bf (6.2.1) Example.} If each $R_i = {\bbC}[Z_i]$ is the  polynomial algebra in 
one 
variable, then $R_1\star...\star R_n = {\bbC}\langle Z_1, ..., Z_n\rangle$ is the 
algebra of noncommutative polynomials. If each $R_i= {\bbC}[X_i, X_i^{-1}]$ us 
the algebra of Laurent polynomials, then $R_1\star ... \star R_n= {\bbC}\langle X_1^{\pm 
1}, 
..., X_n^{\pm 1}\rangle$ is the algebra of noncommutative Laurent polynomials. 

\vskip .1cm 

The following description of the free product follows easily from definition 
(see [VDN]). 

\proclaim (6.2.2) Proposition. Suppose that for each $i$ we chose a subspace 
$R_i^\circ\subset R_i$ which is a complement to ${\bbC}\cdot 1$. Then as a vector 
space 
$$R_1\star ...\star R_n = {\bbC}\cdot 1 \oplus \bigoplus_{k>0}\bigoplus_{i_1\neq i_2\neq 
... 
\neq i_k} R_{i_1}^\circ \otimes ... \otimes R_{i_k}^\circ.$$ 

\vskip .1cm

The following definition of the free product of (noncommutative)  measures is
due to Voiculescu, see [VDN].

\proclaim (6.2.3) Proposition-Definition. Let $R_i$, $i=1, ..., n$ be associative  algebras with 1, 
and $\tau_i: R_i\to \bbC$ be finite normalized measures. Then there exists a unique
finite normalized measure $\tau = \bigstar \tau_i$ on $\bigstar_{i=1}^n R_i$ with the following
properties: \hfill\break
(1) $\tau|_{R_i} = \tau_i.$\hfill\break
(2) If $i_1\neq ... \neq i_k$ and $a_\nu\in R_\nu$ are  such that $\tau_{i_\nu}(a_\nu) = 0$, then
$\tau(a_{i_1} ... a_{i_k}) = 0$. \hfill\break
If the $R_i$ are $*$-algebras, and each $\tau_i$ is a probability measure, then so is $\tau$.

Both the existence and the uniqueness of $\tau$ follow at once from (6.2.2), if we take
$R_i^\circ = \Ker (\tau_i)$. 
The problem of finding $\tau(a_{1} ... a_{k})$ for arbitrary elements $a_{\nu}\in R_{i_\nu}$
is clearly equivalent to that of writing $a_1 ... a_k$ in the normal form (6.2.2). To do this, 
one writes
$$a_\nu = \tau_{i_\nu}(a_\nu) \cdot 1 + a_\nu^\circ,\leqno (6.2.4)$$
with $a_\nu^\circ$ defined so as to satisfy (6.2.4) and we have $\phi_{i_\nu}(a_\nu^\circ) = 0$. Then
one uses the conditions (1) and (2) to distribute. 

\vskip .2cm

\noindent {\bf (6.2.5) Examples.} Suppose we have two algebras $A$ and $B$ and normalized measures
$\phi: A\to \bbC$ and $\psi: B\to\bbC$. Let $\chi: A\star B\to \bbC$ be the free product of $\phi$ and $\psi$.
Then for $a, a'\in A$ and $b, b'\in B$ we have, after some calculations:
$$\chi(ab) = \phi(a)\psi(b), \quad \chi(aba') = \phi(aa')\psi(b),$$
$$\chi(aba'b') = \phi(aa')\psi(b)\psi(b') + \phi(a)\phi(a') \psi(bb') - \phi(a)\phi(a')\psi(b)\psi(b'). $$
See [NS], Thm. 14.4,   for a general formula for $\chi(a_1b_1 ... a_mb_m)$, $a_i\in A$, $b_i\in B$.

\vskip .2cm

\noindent {\bf (6.2.6) Examples.} (a) Let $R_i = \bbC [x^{\pm 1}]$, $i=1, ..., n$
and let $\tau_i$ be given by the integration over the normalized Haar measure $d^*x$ on the unit circle. 
Thus
$$\tau_i \left(f(x) = \sum_m a_m x^m\right)  = \int_{|x|=1} f(x) d^*x = a_0.$$
The free product of these measures is the functional on $\bbC\langle X_1^{\pm1} , ..., X_n^{\pm 1}\rangle$\
given by
$$\tau \left( f(X) = \sum_{\gamma\in F_n} a_\gamma X^\gamma\right) = a_0,$$
the constant term of a noncommutative Laurent polynomial. The asymptotic freedom theorem for unitary matrices (1.4.7) 
says that this functional is the limit, as $N\to\infty$, of the functionals from Example 6.1.1(c)
 with $\mu_N$, for each $N$, 
being the normalized Haar measure on $U(N)^n$. 

(b) Let $R_i = \bbC[z]$, $i=1, ..., n$, and let $\tau_i = 
\delta (z- {a_i})$ be the Dirac delta-function situated at a point $a_i\in \bbC$,
i.e., $\tau_i(f) = f(a_i)$. Then the free product $\tau = \tau_1\star ...\star \tau_n$ is given by:
$$\tau(f(Z_1, ..., Z_n)) = f(a_1\cdot 1 , ..., a_n\cdot 1),$$
in other words, it depends only on the image of $f$ in the ring of commutative polynomials $\bbC[z_1, ..., z_n]$.
This can be seen from the procedure (6.2.4) using the fact that each $\tau_i:  \bbC[z] \to\bbC$ is a ring
homomorphism. 

\vskip .1cm

(c) Let $R_i = \bbC [z]$, $i=1, ..., n$, and let $\tau_i$ be the integration over the standard Gaussian 
probability measure
$$\tau_i(f) = {1\over\sqrt{2\pi}} \int_{-\infty}^\infty f(z) e^{-z^2/2} dz.$$
Their free product is a probability measure on $\bbC\langle Z_1, ..., Z_n\rangle$
denoted by $\xi_n$ and called the {\em free Gaussian measure}. The asymptotic freedom for Hermitian
Gaussian ensembles [V] can be formulated as follows.

\vskip .1cm

\proclaim (6.2.7) Theorem. The measure $\xi_n$ is the limit, as $N\to\infty$,  of the measures from Example 6.1.1(b)
where, for each $N$, we take for $\mu_N $ the   Gaussian  probability measure on the vector space
$(\Herm_N)^n $ corresponding to the
scalar product $\sum \Tr(A_iB_i)$ on this vector space:
$$\mu=\mu_N = {1\over (2\pi)^{n N^2/2}} \exp\left( -{1\over 2} \sum_{i=1}^n Z_i^2\right) dZ_1 ... dZ_n.$$

\vskip .3cm

\noindent {\bf (6.3) The Fourier transform of noncommutative measures.} In classical analysis, Fourier transform is
defined for measures on $\bbR^n$, not on an arbitrary curved manifold. We will call a measure on $\bbR^n_{\NC}$
(``noncommutative $\bbR^n$") a datum consisting of a $*$-algebra $R$, a $*$-homomorphism
$\bbC\langle Z_1, ..., Z_n\rangle \to R$ (i.e., a choice of self-adjoint elements in $R$ which we will still denote $Z_i$),
and a measure $\tau$ on $R$.   Elements of $R$  for which $\tau$ is defined, will be thought of as functions integrable 
with respect to the measure. This concept is thus very similar to that of $n$ noncommutative random variables
in noncommutative probability theory, except we do not require any positivity or normalization.

Let $\tau$ be a measure on $\bbR^n_{\NC}$. Its Fourier transform  is the
 complex valued function $\GF(\tau)$ 
on the group $\Pi_n$
 of piecewise smooth paths in $\bbR^n$ defined as follows: 
$$\GF (\tau) (\gamma) = \tau (E_\gamma(iZ_1, ..., iZ_n)), \quad \gamma\in\Pi_n.\leqno (6.3.1)$$
Here we assume that the ``entire function" $E_\gamma(iZ_1, ..., iZ_n)$ lies in the
domain of definition of $\tau$. In physical terminology, $\GF(\tau)(\gamma)$ is the ``Wilson loop
functional" (defined here for non-closed paths as well). 

\vskip .2cm

\noindent {\bf (6.3.2) Example: delta-functions.}  (a) For every $J= (j_1, ..., j_m)$ we have the
measure $\delta^{(J)}$ on $\bbC\langle \langle Z_1, ..., Z_n\rangle\rangle$ given by
$$\delta^{(J)} \left(\sum_m \sum_{I=(i_1, ..., i_m) } a_I Z_{i_1} ... Z_{i_m}\right) = a_J. $$
The Fourier transform of $\delta^{(J)}$ is the function $W_J: \Pi_n\to\bbC$ which associates to a path $\gamma$
the iterated integral along $\gamma$ labelled by $J$: 
$$W_J(\gamma) = \int^{\rightarrow}_\gamma dy_{j_1} ... dy_{j_m}.$$
We will call these functions {\em monomial} functions on $\Pi_n$

\vskip .2cm

(b) If we take for $\tau$ the free product of (underived) delta-functions $\delta_{a_1} \star ... \star \delta_{a_n}$,
as in Example 6.2.6(b), then 
$$\GF(\tau)(\gamma) = \exp\bigl (i  (e(\gamma), a) \bigr),$$
where $e(\gamma)\in \bbR^n$ is the endpoint of $\gamma$. This follows from the fact that $\tau$
is supported on the commutative locus, i.e., $\tau(E_\gamma(iZ))$ depends only on the image of $E_\gamma(iZ)$ 
in the commutative power series ring, which is
 $\exp\bigl( i  e(\gamma), z\bigr))$. 

\vskip .3cm

\noindent {\bf (6.4) Convolution and product.} Let $\tau,\sigma$  be two measures on $\bbR^n_{\NC}$,
so we have homomorphisms
$$\alpha: \bbC\langle Z_1, ..., Z_n\rangle \to R, \quad \beta: \bbC\langle Z_1, ..., Z_n\rangle \to S$$
and $\tau$ is a linear functional on $R$ while $\sigma$ is a linear functional on $S$. 
Their (tensor) convolution is the measure $\tau *\sigma$ which corresponds to the homomorphism
$$\bbC\langle Z_1, ..., Z_n\rangle \to R\otimes S, \quad Z_i\mapsto \alpha(Z_i)\otimes 1 + 1\otimes\beta(Z_i), \leqno (6.4.1)$$
and the linear functional 
$$\tau * \sigma: R\otimes S\to \bbC, \quad r\otimes s\mapsto \tau(r)\otimes \sigma(s).$$
For commutative algebras this corresponds to the usual  convolution of measures with respect to the group
structure on $\bbR^n$. 

\proclaim (6.4.2) Proposition. The Fourier transform of the convolution of measures is the product of their
Fourier transforms:
$$\GF(\tau * \sigma) = \GF(\tau) \cdot \GF(\sigma).$$

\noindent {\sl Proof:} This is a consequence of the fact that the elements $E_\gamma(iZ_1, ..., iZ_n)$
of $\bbC\langle\langle Z_1, ..., Z_n\rangle\rangle$ are group-like, 
see the exponential property (2.4.7).  \qed 

\vskip .3cm

\noindent {\bf (6.5) Formal Fourier transform of noncommutative measures.} 
The product of two monomial functions on $\Pi_n$ is a linear combination of
monomial functions. This expresses  Chen's shuffle relations among iterated integrals:
$$W_{j_1, ..., j_m} W_{j_{m+1}, ..., j_{m+p}} = \sum_s W_{j_{s(1)}, ..., j_{s(m+p)}},\leqno (6.5.1) $$
the sum being over the set of $(m,p)$-shuffles. An indentical  formula holds for the convolution
of the measures $\delta^{(j_1, ..., j_m)}$  and $\delta^{(j_{m+1}, ..., j_{m+p})}$, as
both formulas describe the Hopf algebra structure on $\bbC\langle\langle Z_1, ..., Z_n\rangle\rangle$. 

The $\bbC$-algebra with basis $W_J = W_{j_1, ..., j_m}$ and multiplication law (6.5.1)
is nothing but the algebra
$$\bbC[G_n] = \lim_{\longrightarrow} \bbC[G_{n,d}] \leqno (6.5.2)$$
of regular functions on the group scheme $G_n = \lim\limits_{\longleftarrow}\,\, G_{n,d}$.
The multiplication in $G_n$ corresponds to the Hopf algebra structure given by
$$\Delta(W_{j_1, ..., j_m}) = \sum_{\nu=0}^{m+1} W_{j_1 ...  j_\nu}\otimes W_{j_{\nu+1}, ..., j_m}.
\leqno (6.5.3)$$
 Elements of $\bbC [G_n]$ can be called polynomial functions on $\Pi_n$.

Note that formal infinite linear combinations (series) $\sum_J c_J W_J$ still form a well defined algebra
via (6.5.1), which we denote $\bbC [[ G_n]]$. This is the algebra of functions on the formal completion
of $G_n$ at 1. The rule (6.5.3) makes $\bbC [[G_n]]$ into a topological Hopf algebra. 

Let $\tau$ be a measure on $\bbR^n_{\NC}$. We will call the formal Fourier transform of $\tau$
 the series
$$\widehat{\GF}(\tau) = \sum_{J=(j_1, ..., j_m)} \tau(Z_{i_1} ... Z_{j_m}) \cdot W_J \quad\in\quad
\bbC [[G_n]]. \leqno (6.5.4)$$
As before, we see that convolution of measures is taken into the product in $\bbC [[G_n]]$.

\vfill\eject

\centerline{\bf  7. Towards the inverse NCFT}

\vskip 1cm

\noindent {\bf (7.0)} In this section we sketch a possible approach to
 the problem of finding the inverse to the NCFT $\cF$ from (2.2).
In other words, given a ``noncommutative function" $f=f(Z_1, ..., Z_n)$, how to find a measure $\mu$
on (possibly some completion of) $\Pi_n$ such that $\cF(\mu)=f$? Note that unlike in classical
analysis, the dual Fourier transform $\GF$ (from noncommutative measures to functions on $\Pi_n$)
does not provide even a conjectural answer, since there is no natural identification of functions
and measures. 

So we take as our starting point the case of discrete NCFT (1.4.5) where Theorem 1.4.6 
provides a neat inversion formula. 

\vskip .3cm

\noindent {\bf (7.1) Fourier series and Fourier integrals.} We recall the classical
procedure expressing  Fourier integrals as scaling limits of Fourier series, see
[W], \S 5. Let $f(x)$ be a piecewise continuous, $\bbC$-valued function on $\bbR$
of sufficiently rapid decay.
We can restrict $f$ to the interval $[-\pi, \pi]$ which is a fundamental domain
for the exponential map $z\mapsto \exp ({iz}), \quad \bbR\to S^1$, and then 
represent $f$ on this interval as  a Fourier series in $e^{inz}$, $n\in \bbZ$. 
 
Next, let us scale the interval to  $[-A, A]$ instead. Then the orthonormal
 basis of functions is 
formed by
$${1\over\sqrt{2A}} \exp\left( {n\pi i  z\over A}\right), 
\quad n\in \bbZ,\leqno (7.1.1)$$
so on the new interval we have
$$f(z) = {1\over 2A} \sum_{n\in\bbZ} \exp \left({n\pi i z\over A}\right) \int_{-A}^A f(w) 
\exp\left( {-n \pi i  w\over A}\right) dy.\leqno (7.1.2)$$
If we associate  the Fourier coefficients to the scaled lattice points, putting
$$g\left( {n\pi\over A}\right)= {1\over\sqrt{2\pi}} \int_{-A}^A f(z)
 \exp\left( {-n\pi i z\over A}\right) dz,\leqno (7.1.3)$$
then 
$$f(z) = {1\over\sqrt{2\pi}} \sum_{n\in \bbZ} g\left( {n\pi\over A}\right) \exp\left({n\pi i z\over A}
\right) \Delta\left({n\pi\over A}\right), \quad z\in [-A, A],\leqno (7.1.4) $$
where $\Delta\left({n\pi\over A}\right) =  {\pi\over A}$ is the step of the dual lattice. So when $A\to\infty$, 
the formulas (7.1.3-4) ``tend to" the formulas for two mutually inverse Fourier transforms
for functions on $\bbR$. In other words, the measures on $\bbR$ (with
 coordinate $y$) given by
infinite combinations of shifted Dirac delta functions: 
$$ {1\over\sqrt{2 \pi}} \,\, {\pi\over A} \sum_{n\in \bbZ}
 g\left( {n\pi\over A}\right) \delta\left(y- {n\pi\over A}\right),
\leqno (7.1.5)$$
converge, as $A\to\infty$, to a measure whose Fourier transform is $f$. 

\vskip .3cm

\noindent {\bf (7.2) Matrix fundamental domains.} We now consider the analog
of the above formalism for Hermitian matrices instead of elements of  $\bbR$ and unitary
matrices instead of those of $S^1$. 
Let $\Herm_N^{\leq A}$ be the set of Hermitian $N$ by $N$
matrices whose eigenvalues  all lie in $[-A, A]$. Then $\Herm_N^{\leq\pi}$ is a fundamental domain for
the exponential map
$$Z\mapsto X= \exp(i Z), \quad \Herm_N\to U(N).\leqno (7.2.1)$$
Note that the Jacobian of the map (7.2.1) is given by
$$J(Z) = \det {} _{N^2\times N^2} {e^{\ad (Z)}-1\over Z } = \prod_{j,k} {e^{i(\lambda_j-\lambda_k)} -1\over \lambda_j-\lambda_k} = 
\prod_{j<k} 2{ 1-\cos(\lambda_j-\lambda_k)  \over (\lambda_k-\lambda_k)^2}. \leqno (7.2.2)$$
Here $\lambda_1, ..., \lambda_N$ are the eigenvalues of $Z$, 
see [Hel], p. 255. Using the formula for the volume of $U(N)$, see, e.g.,  [Mac], 
we can write the normalized Haar measure on $U(N)$ transferred into
$\Herm_N^{\leq\pi}$, as
$$d^*X = { J(Z)dZ \over V_N}, \quad V_N = \prod_{m=0}^{N-1}{2\pi^{m+1}\over m!}. \leqno (7.2.3)$$

Let $f(Z_1, ..., Z_n)$ be a ``good" noncommutative function (for example an entire function
or a rational function defined for all Hermitian $Z_1, ..., Z_n$ and having good decay at infinity).
 Then we can restrict $f$ to $(\Herm_N^{\leq \pi})^n$
and transfer it, via the map (7.2.1),  to a matrix function on $U(N)^n$. This matrix function is clearly nothing
but
$$f(-i \log(X_1), ..., -i \log (X_n)),\leqno (7.2.4) $$
where $-i\log: U(N)\to\Herm_N^{\leq \pi}$ is the branch of the logarithm defined using our choice of
the fundamental domain. Although (7.2.4) is far from being a noncommutative Laurent polynomial 
(indeed, it is typically discontinuous), one can hope to use the procedure of Theorem 1.4.6 to
expand it into a noncommutative {\em Fourier series}. In other words, assuming that for each
$\gamma\in F_n$ the limit
 $$a_\gamma =  \lim_{N\to\infty} {1\over N} 
 \Tr\int_{X_1, ..., X_n\in U(N)} f(-i\log(X_1), ..., -i\log(X_n)) X^{-\gamma} \prod_{j=1}^n d^*X_j = 
 \leqno (7.2.5)$$
$$ = \lim_{N\to\infty} {1\over N} \Tr\int_{Z_1, ..., Z_n\in\Herm_N^{\leq\pi}} 
f(Z_1, ..., Z_n) E_{\gamma^{-1}}(iZ_1, ..., iZ_n))\prod_{j=1}^n{ J(Z_j) dZ_j\over V_N}$$
exists, we can form the series
$$\sum_\gamma a_\gamma X^\gamma \quad =\quad \sum_\gamma a_\gamma \,\, E_\gamma(iZ_1, ..., iZ_n), \quad Z_j\in \Herm_N^{\leq \pi}.
\leqno (7.2.6)$$
By analogy  with  the classical case one can expect  that this series converges to $f|_{(\Herm_N^{\leq\pi})^n}$ away from the
boundary.

\vskip .3cm

\noindent {\bf (7.3) Scaling the period.} 
 In the situation of (7.2) let us choose $A>0$  and  restrict $f$ to $(\Herm_N^{\leq A})^n$. 
 The same procedure would then expand the restriction into
 a series in
 $$X_j^{\pi/A} = \exp\left(i {\pi\over A} Z_j\right), \quad j=1, ..., n. \leqno (7.3.1)$$
 Let $F_n^{\pi/A}\subset G_n(\bbR)$ be the group generated by the $X_j^{\pi/A}$.
 We can think of elements of $F_n^{\pi/A}$ as rectangular paths in $\bbR^n$ with increments being integer
 multiples of $\pi/A$. The coefficients
 of the series for the restriction given then a function
 $$g_A: F_n^{\pi/A}\to \bbC,$$
so the series will have the form
$$f(Z) = \sum_{\gamma\in F_n^{\pi/A}} g_A(\gamma) E_\gamma(iZ), \quad Z=(Z_1, ..., Z_n), \,\, Z_j\in\Herm_N^{\leq A}. $$
Now, as $A\to\infty$, we would like to say that the $g_A$, considered as linear combinations of Dirac measures
on $\Pi_n$ (or some completion) tend to a limit measure. Although $\Pi_n$ is not a manifold, we can pass to finite
dimensional approximations
$$\Pi_n \subset G_n(\bbR)\buildrel p_d\over\longrightarrow  G_{n,d}(\bbR).$$
Let $F_{n,d}^{\pi/A} = p_d(F_n^{\pi/A})$. This is a free nilpotent group of degree $d$ on generators
$p_d(X_j^{\pi/A})$, and is a discrete subgroup (``lattice") in $G_{n,d}(\bbR)$. As $A\to\infty$, these lattices
are getting dence in $G_{n,d}(\bbR)$. 
Supposing
that the direct image (summation over the fibers) $p_{d*}(g_A)$ exists as a function on $F_{n,d}^{\pi/A}$ or,
what is the same, a
 measure on $G_{n,d}(\bbR)$ supported on the discrete subgroup $F_{n,d}^{\pi/A}$, we can then ask for the existence of
 the limit
 $$\mu_d = \lim_{A\to\infty} p_{d*}(g_A) \quad\in\quad {\operatorname{Meas}}(G_{n,d}(\bbR)).$$
 These measures, if they exist, would then form a pro-measure $\mu_\bullet$  which is the natural
 candidate for the inverse Fourier transform of $f$. 
 The author hopes to address these issues in a future paper.

\vfill\eject

\centerline{\bf REFERENCES}

\vskip 1cm

\noindent [Ba] F. Baudoin, {\it An Introduction to the Geometry of Stochastic Flows},
World Scientific, Singapore, 2004. 

\vskip .1cm

\noindent [Bel] D. R.  Bell, {\em 
Degenerate Stochastic Differential Equations and Hypoellipticity,  }
(Pitman Monographs and Surveys in Pure and Applied Mathematics, 79)  Longman, Harlow, 1995.

\vskip .1cm

\noindent [Be] G. Ben Arous, {\em Flots et s\'eries de Taylor stochastiques}, 
J. of Prob. Theory and Related Fields,  81 (1989), 29-77. 

\vskip .1cm

\noindent [BGG] R. Beals, B. Gaveau, P.C.  Greiner, 
{\em Hamilton-Jacobi theory and the heat kernel on Heisenberg groups},
J. Math. Pures Appl. 79 (2000), 633-689. 

\vskip .1cm

\noindent [Bi1] J.-M. Bismut, {\em Large Deviations and the Malliavin Calculus},
Birkhauser, Boston, 1984.

\vskip .1cm

\noindent [Bi2] J.-M. Bismut, {\em The Atiyah-Singer theorem: a probabilistic approach. I. The index theorem},
  J. Funct. Anal.  57  (1984),   56-99.

\vskip .1cm

\noindent [C0] K.-T.-Chen, {\em Collected Papers},  Birkhauser, Boston, 2000. 

\vskip .1cm

\noindent [C1] K.-T. Chen, {\em Integration in free groups}, Ann. Math. 
54 (1951), 147-162.

\vskip .1cm

\noindent [C2] K.-T. Chen, {\em Integration of paths--a faithful representation
of paths by noncommutative formal power series}, Trans. AMS, 89 (1958), 395-407. 

\vskip .1cm

\noindent [C3] K.-T. Chen, {\em Algebraic paths}, J. of Algebra, 10 (1968), 8-36. 

\vskip .1cm

\noindent [Coh] P. M. Cohn,   {\em Skew fields. Theory of general division rings}, 
 (Encyclopedia of Mathematics and its Applications, 57),   Cambridge University Press, Cambridge, 1995.
 
 \vskip .1cm
 
 \noindent [Con] A. Connes, {\em Noncommutative Geometry},  Academic Press, Inc. San Diego, CA, 1994.

\vskip .1cm

\noindent [FN] M. Fliess, D. Normand-Cyrot, {\em Alg\`ebres de Lie nilpotents, formule de Baker-Campbell-Hausdorff et
int\'egrales iter\'ees de K. -T. Chen}, S\'em. de Probabilit\'es XVI, Springer Lecture
Notes in Math. 920, p. 257-267. 

\vskip .1cm

\noindent [Fo] G. B.  Folland, {\em A fundamental solution for a subelliptic operator}, 
Bull. AMS, 79 (1973), 367-372. 

\vskip .1cm

\noindent [G] B. Gaveau, {\em Principe de moindre action, propagation de chaleur
et estim\'ees sous-elliptiques pour certains groupes nilpotents},
Acta Math. 139 (1977), 95-153. 

\vskip .1cm

\noindent [HL]  M. Hakim-Dowek, D. L\'epingle, {\em L'exponentielle stochastique des
groupes de Lie}, S\'em. de Probabilites XX, Springer Lecture Notes in Math. 1203,
p. 352-374. 

\vskip .1cm

\noindent [Hel] S. Helgason, {\em Differential Geometry, Lie Groups and Symmetric Spaces,}
Academic Press, 1962. 

\vskip .1cm

\noindent [HP] F. Hiai, D. Petz, {\em The Semicircle Law, Free Random Variables and
Entropy}, Amer. Math. Soc. 2000. 

\vskip .1cm

\noindent [Ho] L. H\"ormander, {\em Hypoelliptic second order differential
equations,} Acta Math. 119 (1967), 141-171. 

\vskip .1cm

\noindent [Hu] A. Hulanicki, {\em The distribution of energy in the Brownian
motion in the Gaussian field and analytic hypoellipticity of certain
subelliptic operators on the Heisenberg group}, Studia Math. 
56 (1976), 165-173. 

\vskip .1cm

\noindent [K1] M. Kapranov, {\em Free Lie algebroids and path spaces}, 
in preparation.

\vskip .1cm

\noindent [K2] M. Kapranov, {\em Membranes and higher groupoids},
in preparation. 

\vskip .1cm

\noindent [KW] N. Kunita, S. Watanabe, {\em Stochastic Differential Equations and Diffusion 
processes}, Noorth-Holland, Amsterdam, 1989. 

\vskip .1cm

\noindent [Mac] I. G. Macdonald, {\em The volume of a compact Lie group}, Invent. Math. 56 (1980),
93-95. 

\vskip .1cm

\noindent [McK] H. P. McKean, {\em Stochastic Integrals,} Chelsea Publ. Co. 2005.   

\vskip .1cm

\noindent [Ma] Y. I.  Manin, {\em Iterated Shimura integrals}, preprint
math.NT/0507438. 

\vskip .1cm

\noindent [NS] A. Nica, R. Speicher, {\em Lectures on the Combinatorics of Free Probability},
(London Math. Soc. Lecture Notes vol. 335), Cambridge Univ. Press, 2006. 

\vskip .1cm

\noindent [Ok] B. Oksendal, {\em Stochastic Differential Equations}, Springer, Berlin, 1989. 

\vskip .1cm

\noindent [Pa] A. N. Parshin, {\em On a certain generalization of the Jacobian manifold},
Izv. AN SSSR, 30(1966), 175-182.

\vskip .1cm

\noindent [Po] A.M. Polyakov, {\em Gauge Fields and Strings}, Harwood Academic Publ. 1987.

\vskip .1cm

\noindent [R] C. Reutenauer, {\em Free Lie Algebras}, Oxford Univ. Press, 1993. 

\vskip .1cm

\noindent [Si] I. M. Singer, {\em On the master field in two dimensions,} in: Functional Analysis
on the Eve of the 21st Century (In honor of I. M. Gelfand), S. Gindikin et al. Eds, Vol.1, 263-281.  

\vskip .1cm

\noindent [SW] D. W. Stroock, S. R. S. Varadhan, {\em Multidimensional Diffusion Processes},
Springer, Berlin, 1979.  

\noindent [Ta] J. L. Taylor, {\em Functions of several noncommuting variables,} Bull. AMS,
79 (1973), 1-34. 

\vskip .1cm

\noindent [V] D. Voiculescu, {\it Limit laws for random matrices and free products},
Invent. Math. 104 (1991), 201-220. 

\vskip .1cm

\noindent [VDN] D. Voiculescu, K.J. Dykema, A. Nica, {\it Free Random Variables},
Amer. Math. Soc. 1992.  

\vskip .1cm

\noindent [W] N. Wiener, {\em The Fourier Integral and Certain of Its Applications,}
Dover Publ, 1958.

\vskip 1cm

\end{document}